\documentclass[12pt]{article}
\usepackage{authblk,graphicx,chemarr,amsmath,amssymb,amsthm,mathrsfs}
  \graphicspath{{fig/}}
 \renewcommand\det{\mbox{det}}
  \usepackage[margin=1in]{geometry}
  \newcommand\ab[1]{\vert{#1}\vert}
    \newcommand\G{\Gamma}
 \newcommand{\rank}{\mathop{\mbox{rank}}}

\newcommand{\lra}{\mathop{\longrightarrow}\limits}
\newcommand{\Rnn}{{\mathbb R}_{\ge 0}^n}

\newcommand{\sgn}{\mathop{\mbox{sgn}}}

\newcommand{\tV}{{\tilde{V}}}

\newcommand{\cW}{{\mathcal W}}

\newcommand{\R}{{\rm \bf R}}
\newcommand{\Rs}{{\mathscr R}}
\newcommand{\Ss}{{\mathscr S}}
\newcommand{\Ns}{{\mathscr N}}
\newcommand{\tl}{\tilde}
\newcommand\comment[1]{}

 \theoremstyle{thmstyleone}%
 \newtheorem{theorem}{Theorem}%
  \newtheorem{lemma}[theorem]{Lemma}%
  \newtheorem{corollary}[theorem]{Corollary}%

 \theoremstyle{thmstyletwo}%
 \newtheorem{remark}{Remark}%
 
 \theoremstyle{thmstylethree}%
 \newtheorem{definition}{Definition}%
 
\begin{document}
 
 \title{\Huge Graphical characterizations of robust stability in biological interaction networks}

 	\author[1]{{M.} {Ali Al-Radhawi}} 
 	
 	\affil[1]{{Department of Electrical \& Computer Engineering, and the Center of Theoretical Biological Physics},  {Northeastern University},  { {360 Huntington Avenue},  {Boston}, {02115},  {Massachusetts}, {United States}}.~Email:~\texttt{{malirdwi@northeastern.edu}} }

 \date{November 20th, 2022}
 
 \maketitle

	 \begin{abstract}{	Previous studies have inferred robust stability of reaction networks by utilizing linear programs or iterative algorithms. Such algorithms become tedious or computationally infeasible for large networks. In addition, they operate like black-boxes without offering intuition for the structures that are necessary to maintain stability. In this work, we provide several graphical criteria for constructing robust stability certificates, checking robust non-degeneracy, {verifying  persistence},  and establishing global stability. By characterizing a set of stability-preserving graph modifications that includes the enzymatic modification motif, we show that the stability of arbitrarily large nonlinear networks can be examined by simple visual inspection. We show  applications of this technique to ubiquitous motifs in systems biology such as  Post-Translational Modification (PTM) cycles, the Ribosome Flow Model (RFM), \emph{T}-cell kinetic proofreading and others. The results of this paper are dedicated in honor of Eduardo D. Sontag's  seventieth birthday and his pioneering work in nonlinear dynamical systems and  mathematical systems biology.}\end{abstract}

	
	\section{Introduction}
	Biomolecular Interaction Networks (BINs) function under severe forms of external and internal uncertainty. Nevertheless, they operate robustly and consistently to maintain \emph{homeostasis}, which is understood as the maintenance of a desired steady-state against environmental factors, external signals, and \emph{in-vivo} fluctuations in the concentrations of   biochemical species. In fact, robustness has been proposed as a key \emph{defining} property of biological networks \cite{morohashi02,kitano02}. However, mathematical analysis of such networks has been lagging as the dynamical system descriptions of such networks suffer from  \emph{nonlinearity} and \emph{uncertainty}. {Generic} nonlinear dynamical systems are already difficult to analyze due to {the scarcity of general and powerful analysis tools. Furthermore, they can manifest complex forms of unstable behavior that are not exhibited by linear systems. For instance,}    small fluctuations in   concentrations, or tiny changes in kinetic parameters, can have radical effects causing the observable phenotype to be driven to a different region of the state space, and/or to lose stability altogether and transform into a sustained oscillation or chaotic behavior. This may make the biological network lose its function and cause key species to reach undesirable or even  {unsafe} levels. In fact, disease can be often characterized mathematically as the loss of stability of a certain phenotype \cite{maclean15,langlois17}. A second complicating factor is the fact that the exact form of kinetics (determining the speed of interactions) are difficult to measure and are subject to environmental changes.
	\begin{figure}
		\centering
		\includegraphics[width=0.67\columnwidth]{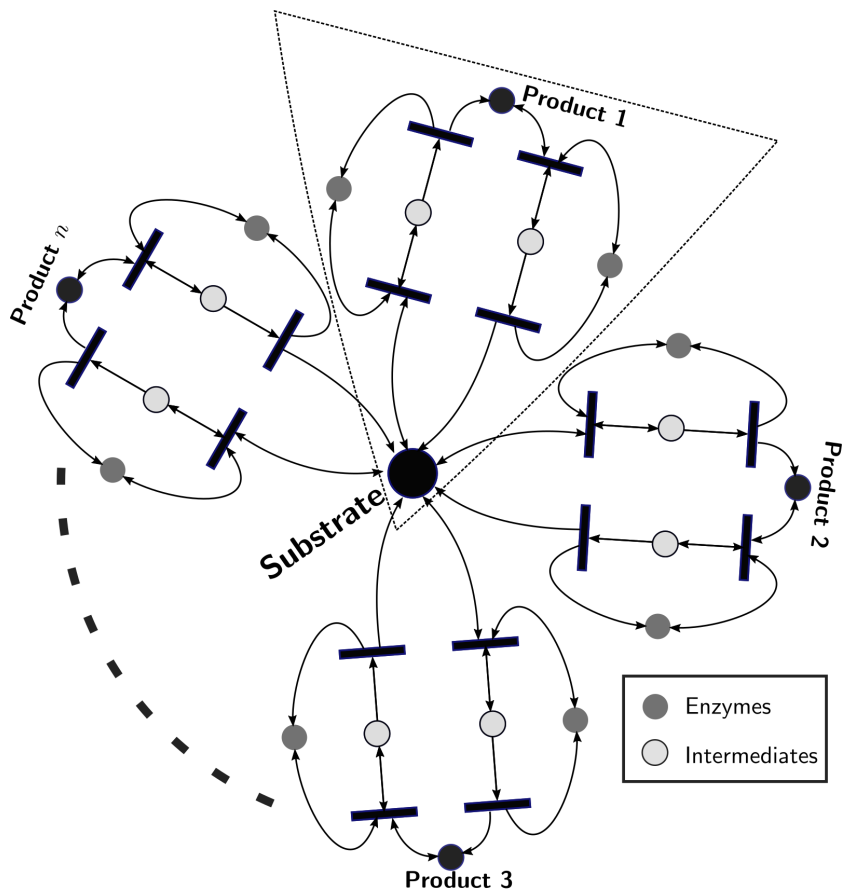}
		\caption{\textbf{Computational RLF construction is tedious for large networks.} The figure depicts a Petri-net representation of the \textit{PTM star}: a substrate that is a target of an arbitrary finite number of  distinct competing PTM    cycles (e.g, phosphorylation, methylation, ubiquitination, etc). The subnetwork inside  the dotted triangle depicts a single PTM cycle. A rectangle denotes a reaction, while a circle denotes a species.   }
		\label{f.star}
	\end{figure}
	Therefore, verifying the stability of a given nonlinear BIN without reference to its kinetics has been a challenging long-standing goal  in systems biology research \cite{bailey01}. Nevertheless, partial success has been achieved in this endeavor.  Examples  include the theory of complex balance \cite{horn72}, \cite{feinberg87}, \cite{sontag01}, and the theory of monotone BINs \cite{angeli10}. More recently, stability certificates have been constructed via  Robust Lyapunov Functions (RLFs) in reaction \cite{MA_cdc13,PWLRj,MA_LEARN}, and concentration coordinates \cite{MA_cdc14,blanchini14,blanchini17,MA_LEARN,blanchini21}. Except for a small subclass of BINs (see  §III.C), such methods mainly utilize computational algorithms to construct   RLFs via either iterative algorithms or linear programs.
	However, such algorithms act as   “black-boxes” and are not interpretable in terms of the structural properties of the network's graph. This has several drawbacks.
	First,   computational algorithms become tedious for larger networks as the number of species  and reactions grow. Consider the \emph{PTM star} depicted in Figure \ref{f.star} whose size grows considerably for large $n$. 
	Second,  ``stability-preserving'' graph modifications are not well characterized. A simple modification of the BIN graph mandates a re-run of the computational algorithm from scratch. For instance, is the stability of the PTM star preserved if we added  inflow/outflow reactions for the substrate ($\emptyset \rightleftharpoons \mbox{Substrate}$)?. Third, fundamental ``motifs'' have been described as the building blocks of BINs \cite{alon06}. However, a corresponding ``modular'' theory for RLF construction that utilizes the stability properties of its subnetworks is lacking.  For example, the difference between the PTM star (Fig. \ref{f.star}) with $n$ and $n+1$ products  is in the addition of an extra \emph{PTM cycle}. How does the addition of the extra motif affect stability? 
	
	The above questions are hard to answer using computational algorithms. In this work, we identify a set of   stability preserving graph modifications. In particular, we show  that  the stability of many large networks in  systems biology     can be understood modularly. For the specific network in Figure \ref{f.star}, we will show that it can be ``reduced'' to a simple linear network (See Figure \ref{fig:ptmstarlinear} {in sec. \ref{s.ptmstar}}). Hence, it admits a stability certificate for every $n \ge 1$, a result which is not readily achievable using previous results \cite{horn72,feinberg87,sontag01,angeli10,blanchini14,MA_LEARN}. We will show that   the addition of an inflow/outflow reaction to the substrate preserves stability, and that the PTM cycle is a fundamental ``stable'' motif in a precise manner to be defined.
	
	Our unified framework can be applied to many networks in the literature whose stability was studied individually via various techniques, this includes the T-cell kinetic proofreading network \cite{sontag01}, the PTM cycle \cite{angeli08}, the all-encompassing processive PTM cycle \cite{shiu17},  the ribosome flow model and its variations \cite{rfm12,rfm16,LCSS_jared}, and others. 
	
	It is worth noting that many of the properties of BINs have already been characterized graphically. This includes complex balance \cite{horn72},\cite{feinberg87}, injectivity \cite{craciun05,banaji16},  monotonicity \cite{angeli10,defreitas17}, and persistence \cite{angeli07p}. {Additional studies have tackled graph modifications that preserve various other properties of BINs \cite{gross20,banaji22}.} Therefore, we complement this literature by characterizing robust stability in graphical terms for classes of BINs {for the first time.}
	
	The paper proceeds as follows. Section II reviews notation and definitions. Section III reviews relevant results on linear (mono-molecular) networks. In section III, we list the graph modifications under consideration, and show the existence of RLFs for classes of modified networks. Global stability and robust non-degeneracy are discussed in section IV. Applications are studied in section V. Proofs are included in the appendix.
	\section{Background and Notation}
	\subsection{Biological Interaction Networks}
	Any collection of chemical reactions can be written mathematically using the formalism of of Biological Interaction Networks (BINs). Hence, we review the standard definitions and notation  \cite{feinberg87,erdi89,sontag01,angeli09tut,MA_LEARN}.
	
	A BIN (also known as  a Chemical Reaction Network (CRN)) is a pair $\Ns=(\Ss,\Rs)$, where  $\mathscr S=\{X_1,..,X_n\}$ is the set of species, and  $\Rs=\{\R_1,...,\R_\nu\}$ is the set of reactions.
	A species is the entity that partakes in  or is formed in a chemical interaction. Within the realm of biomolecular networks, a species can be a substrate, a  complex, an enzyme, an mRNA molecular, a gene promoter state, etc.   A reaction is the transformation of reacting species into product species. Examples include complex formation, binding, unbinding, decay, production, complex formation, etc. 
	
	The mathematical structure of BINs  can be described by two mathematical substructures: \emph{the stoichiometry} and \emph{the kinetics}.

	\subsubsection{The Stoichiometry} The relative gain or loss of molecules of species $X_i$ between the sides of each reaction is the \emph{stoichiometry of $X_i$}. This is represented by writing a reaction as: 
	\begin{equation}\label{e.reaction}
		\R_j: \quad \sum_{i=1}^n \alpha_{ij} X_i \longrightarrow \sum_{i=1}^n \beta_{ij} X_i, \ j=1,..,\nu,
	\end{equation}
	where $\alpha_{ij}, \beta_{ij} \ge 0$ are integers known as the \emph{stoichiometry coefficients}. %
	If a transformation can happen also in the reverse direction, then $\R_j$ is said to be \emph{reversible} and its reverse is denoted by $\R_{-j}$. A reaction can have no reactants or no products (though not simultaneously). The empty side is denoted by $\emptyset$. 
	
	If a reaction has a species both as a reactant and as a product (for example,  $X+Y \to X$) then it is called \emph{catalytic}. %

	The \emph{stoichiometry matrix} $\Gamma$ of a given network is an $n \times \nu$ matrix whose $(i,j)$th entry describes the net gain/loss of the $i$th species at the $j$th reaction. Hence, it can be written element wise as:
	$   [\Gamma] _{ij} = \beta_{ij}- \alpha_{ij}.$

	\subsubsection{Kinetics} The set of relationships  that determine the speed of transformation of reactant species into product species are known as \emph{kinetics}.    In order to describe such relations, the species need to be quantified. 
	A species $X_i$ is quantified by assigning it a non-negative real number known as the concentration $x_i {\in \mathbb R_{\ge0}^n}${, where $\Rnn$ denotes  the non-negative orthant in the $n$-dimensional Euclidean space.} A reaction $\R_j$ is assigned a single-valued mapping $R_j:{\Rnn \to  {\mathbb R}_{\ge 0}}$  known as the \emph{reaction rate}. The reaction rate vector is written as $R(x)=[R_1(x),...,R_\nu(x)]^T$.    %
	
	The most common form of kinetics is known as \emph{Mass-Action} and it can be written as:
	$\label{e.mass_action}
	R_j(x)= k_j \prod_{i=1}^n x_i^{\alpha_{ij}}$, where
	$k_j>0, j=1,..,\nu$ are  the \textit{kinetic constants}. However, this form ``is not based on fundamental laws'' and is merely ``good phenomenology'' justified by imagining the reactants as colliding molecules \cite{gunawardena14}. In biological systems, in particular, other forms of kinetics usually arise when modeling networks involving multiple time-scales. This includes Michaelis-Menten, Hill kinetics, etc.
	Therefore,  we do not assume a specific functional form of kinetics. We only assume that the kinetics are \emph{monotone}. More precisely,    the reaction rates $R_j(x), j=1,..,\nu$  satisfy:
 	\begin{enumerate}
		\item[\!\!] {\bf AK1}. each reaction varies smoothly with respects to its reactants, i.e $R(x)$ is $C^1$;
		\item[\!\!] {\bf AK2}. a reaction requires all its reactants to occur, i.e.,  if $\alpha_{ij}>0$, then $x_i=0$ implies  $R_j(x)=0$; %
		\item[\!\!] {\bf AK3}. if a reactant increases, then the reaction rate increase, i.e ${\partial R_j}/{\partial x_i}(x) \ge 0$ if $\alpha_{ij}>0$ and ${\partial R_j}/{\partial x_i}(x)\equiv 0$ if $\alpha_{ij}=0$. Furthermore, the aforementioned inequality is strict whenever the reactants are strictly positive.
	\end{enumerate}
	For a given network $\Ns$, the set of a reaction rates satisfying the assumptions above is called the \emph{admissible kinetics}. Furthermore,  the assumptions AK1-AK3 translate into a sign-pattern constraint on the Jacobian of $R$. To formalize this, %
	 let  $\mathcal K_\Ns$ be defined as $\mathcal K_\Ns=\{V\in \mathbb R^{\nu \times n} \vert [V]_{ji}>0~\mbox{whenever}~  X_i~\mbox{is a reactant of}~\R_j, ~\mbox{and}~[V]_{ji}=0~\mbox{otherwise}\}$. We think of   $\mathcal K_\Ns $ as the set of all possible Jacobian matrices $\partial R/\partial x$ evaluated on the \textit{positive} orthant {$\mathbb R_+^n$}.    %

	\subsubsection{Dynamics}   We view the concentrations as trajectories in time and write them as $x(t)=[x_1(t),...,x_n(t)]^T$. 
	{The temporal evolution of the network is given by the following Ordinary Differential Equation (ODE)}:
	\begin{equation} 
		\label{e.ode} 
		\dot x = \Gamma R(x), \; \; x(0)=x_\circ.
	\end{equation}
		The positive orthant is forward-invariant for \eqref{e.ode}, i.e. if $x_\circ$ is positive, then the trajectory stays positive for all time $t \ge 0$.

	In the biomolecular context, there are usually \emph{conserved quantities} which do not get created or annihilated during the course of the reaction. This can include total amounts of DNA, enzymes, substrates, ribosomes, etc. 
	Mathematically, a stoichiometric conservation law is a nonnegative vector $d \in \mathbb R_{\ge 0}^n$ satisfying  $d^T \Gamma =0$. If $d$ is positive then the network is called \emph{conservative}.
	
	The existence of a conservation law implies that $d^Tx(t)\equiv d^T x(0)$. Hence, the positive orthant is partitioned into a foliage of subsets known as \emph{stoichiometric classes}. For  a state vector $x_\circ$, the corresponding class is written as $\mathscr C_{x_\circ}:=(\{x_\circ\}+ \mbox{Im}(\Gamma)) \cap \Rnn$, and it is forward invariant. Therefore, all Lyapunov functions and claims of stability are relative to a stoichiometric class. 
	For a conservative network, all stoichiometric classes are compact polyhedral sets, and hence all trajectories are bounded. In addition, this  guarantees   at least one steady state in each  stoichiometric class by applying Brouwer's fixed point theorem to the associated flow of the dynamical system restricted to the stoichiometric class.%
	
	A vector $v$ is called a \textit{flux} if $\Gamma v=0$.  In order to simplify the treatment, we will  assume the following about the stoichiometry of the network:
	\begin{enumerate}
		\item[ ] {\textbf{AS1.}} There exists a positive flux, i.e.,  $\exists v \in \ker \Gamma$ such that $v\gg0$.  
		\item[] {\textbf{AS2.}} The network has no catalytic reactions.
	\end{enumerate}
	Assumption AS1 is necessary for the existence of positive steady states for the corresponding dynamical system \eqref{e.ode}.

	\subsection{Graphical representation: Petri-Nets}
	
	BINs can be represented {graphically} in several ways. We adopt the \emph{Petri-net} formalism \cite{petri08} (also known as the species-reaction graph \cite{craciun06}). A Petri-net is a weighted directed bipartite graph. The vertices consists of the set of   species $\Ss$ (represented by circles) and the set of reactions $\Rs$ (represented by rectangles). An edge with a weight $w$ from $X_i\in\Ss$ to $\R_j\in\Rs$ means that $X_i$ is a reactant of $\R_j$ with stoichiometric coefficient $w$, while the reverse edge means that $X_i$ is a product of $\R_j$ with a stoichiometric coefficient $w$. For a more compact representation, if two reactions are the reverse of each other (e.g, $\R_j,\R_{-j}$) then they are represented as a single reaction in the Petri-net with reversible edges.	 In the formalism of Petri-nets \cite{murata89}, the stoichiometric matrix $\Gamma$ is the \emph{incidence matrix} of the Petri-net. 

	For example, the PTM star in Fig. 1 corresponds to the following network:
	\begin{align} \label{starA}
		S+E_i & \rightleftharpoons C_i \longrightarrow P_i + E_i, \\ \label{starB}
		P_i+F_i & \rightleftharpoons D_i \longrightarrow S + F_i,
	\end{align}
	$i=1,..,n$, where $S$ denotes the substrate and $P_i$ denotes the $i$th product.
	
	\subsection{Robust Lyapunov Functions}
	Following our previous work \cite{MA_cdc14,PWLRj,MA_LEARN}, a locally Lipschitz function $V: \mathbb R^n \to \mathbb R_{\ge 0}$ is a Robust Lyapunov Function (RLF) for a given network $\Ns$  iff:
	\begin{enumerate} \item it is \emph{positive-definite}, i.e., $V(x) \ge 0$ for all $x$, and $V(x)=0$ iff $\Gamma R(x)=0$, and
		\item it is \emph{non-increasing}, i.e., $\dot V(x) \le 0$ for all $x$ and all $R$ {satisfying $\frac{\partial R}{\partial x} \in\mathcal K_\Ns$}.
	\end{enumerate}
	Since $V$ is not assumed to be continuously differentiable, the derivative above is defined in the sense of Dini  as $\dot V(x):=\limsup_{h \to 0^+} ( V(x+h\Gamma R(x))-V(x))/h$ \cite{yoshizawa}. Existence of an RLF guarantees that the steady state set is Lyapunov stable, and   that all $V$'s level sets are trapping \cite{yoshizawa,PWLRj,MA_LEARN}. Global stability can be verified by a LaSalle argument or by establishing robust non-degeneracy of the Jacobian \cite{PWLRj,blanchini17,MA_LEARN}.  In this paper, we utilize RLFs that can be written as piecewise linear (PWL) functions in terms of the rates. In \cite{MA_LEARN}, it has been shown that they can be converted to PWL RLFs in the concentration-coordinates and vice versa. Hence, we will subsequently use the term ``PWL RLF'' to designate an RLF that is piecewise linear either in the rates or in the concentrations.

	\section{Linear (Mono-molecular) Networks} 
	\subsection{Definition {and review}}
	
	Studying general nonlinear BINs is, predictably, a difficult and open problem. In comparison, assuming linearity simplifies the analysis considerably. In order to get a linear ODE  with Mass-Action kinetics, all the reactions have to be \emph{monomolecular}. In other words, there is only a unique reactant with stoichiometry coefficient 1 for each reaction. The resulting ODE can be studied via standard analysis methods for positive linear systems \cite{luenberger79,haddad09}, or as a special case of complex-balanced networks \cite{feinberg87}.    {A weaker notion of linearity is a graphical one where the Petri-net is assumed to be linear, \cite{marinescu91}, which means that each reaction has a unique reactant and a unique product with the stoichiometry coefficients equal to one. Therefore, nonlinear reaction rates are allowed. It has been long-observed that the linearity of the Petri-net is sufficient  for  analysis, i.e., stability analysis can be performed for general monomolecular networks with monotone kinetics \cite{maeda78}. This generalized class of networks is often known as \emph{compartmental networks} \cite{jacquez93}. Hence, we refer to such networks as \emph{linear networks} since the corresponding Petri-net is linear.}  Therefore, we use a \textit{graphical notion of linearity} and not a kinetic one.  The definition is stated formally below:
	
	\begin{definition} A given BIN $\Ns$ is said to be \emph{linear} if each reaction can be written as either $X_i \lra^{R_{ij}} X_j$,  $\emptyset \lra^{u_i} X_i$, or $X_i \lra^{R_{i}} \emptyset$ for some $i,j$ where $R_{ij},R_{i}:\mathbb R_{\ge 0} \to \mathbb R_{\ge 0}, u_i\ge 0$ are the reaction rates.	
	\end{definition} 
	Applying the assumptions AK1-AK3, we note that $R_{ij}$ can be any single-valued strictly increasing $C^1$ function that vanishes at the origin. %
	\subsection{Existence of Lyapunov functions: Sum-of-Currents (SoC) RLF}
	
	One of the advantages of studying linear networks is that their stability is well-characterized. Indeed, it has been long-known 
	\cite{maeda78,jacquez93} that linear networks can be studied using a Lyapunov function of the form: 
	\begin{equation}\label{e.lyap_linear}
		V(x) = \| \dot x \|_1= \sum_{i=1}^n \left \vert \sum_{j \ne i} (R_{ji}(x_j) - R_{ij}(x_i)) + u_i - R_i(x_i) \right \vert,
	\end{equation}
{where $u_i \ge 0$ is the inflow to species $X_i$.}
Note that $V$ is PWL in terms of the rates. 
	We state the following theorem {that restates the result in \cite{maeda78} using our terminology}:
	\begin{theorem}\label{th1}
		Let $\Ns$ be a linear BIN    with any set of admissible reaction rates $\{R_{ij}(x_i),$ $ R_i(x_i),u_i\}_{i,j=1}^n$. Let \eqref{e.ode} be the associated ODE. Let $V$ be defined as in \eqref{e.lyap_linear}. Then,  $V$ is an RLF for $\Ns$.
	\end{theorem}
	
	In order to generalize the result above to classes of nonlinear networks, we will provide a new proof of Theorem 1 in the Appendix based on the techniques used in \cite{PWLRj,MA_LEARN,LCSS_jared}. The same techniques will be generalized to prove Theorem \ref{mainth}. 
	In \cite{MA_LEARN}, we have called \eqref{e.lyap_linear} a Sum-of-Currents (SoC) RLF, since it is a sum of the absolute values of the currents $dx_i/dt, i=1,..,n$, which is analogous to the electric current $I=dq/dt$, where $q$ is the electric charge. 

\subsection{Existence of Lyapunov functions: Max-Min RLF}
For a subclass of linear BINs, another Lyapunov function can be used to   establish stability, which is the Max-Min RLF \cite{MA_cdc13,PWLRj}. Define the set-valued function: $\mathcal R(x)=\{R_{ij}(x),R_i(x),u_i \vert i,j=1,..,n, i\ne j \} $. Then, consider the following function:
\begin{equation}\label{maxmin} V(x)= \max \mathcal R(x) - \min \mathcal R(x),\end{equation}
Note that $V$ is PWL in terms of the rates. 
The existence of an RLF of the form \eqref{maxmin} can be characterized graphically for general BINs \cite{MA_cdc13,PWLRj}. In order to minimize the notational inconvenience, we assume that $\mathbf 1$ is a flux for the network $\Ns$. Hence, the result can be stated as follows:
\begin{theorem}[\cite{MA_cdc13,PWLRj}]\label{th.maxmin} 
Let  a BIN $\Ns$ be given. Assume that it has a \textit{unique} positive flux equal to $\mathbf 1$ and every species $X_i$ is a reactant to a unique reaction. Then, $V$ as defined  in \eqref{maxmin} is an RLF for $\Ns$.
\end{theorem}
\begin{remark}In order to generalize Theorem \ref{th.maxmin} to accommodate BINs that admit a unique positive flux $v\gg 0$, the reactions in $\mathcal R(x)$ can be weighed by the  corresponding entry in $v$ \cite{PWLRj}.\end{remark}
\section{Stability-Preserving Graph Modifications}
\subsection{Definitions}
Consider a BIN $(\mathscr S,\mathscr R)$ that admits an RLF $V$. Assume that the network is modified to a new network $(\tl \Ss, \tl \Rs)$. We are interested in the existence of an RLF for the new network. To be more concrete, we focus on  graph modifications   listed in Table I. As can be noticed, some of these modification can change a linear network into a nonlinear network.  First, we formalize the concept of adding an extra product or reactant to a reaction.
\begin{definition}
Consider a BIN $(\Ss,\Rs)$. We say that a reaction $\tl\R_j$ is an \emph{extension} of a reaction $\R_j \in \Rs$ if   the following holds for each $ X_i \in \Ss$: if  $X_i \in \Ss$ is a reactant  of $\R_j \in \Rs$, then $X_i$ is a reactant of $\tl\R_j \in \tl\Rs$ with the same stoichiometric coefficient. Similarly,  if $X_i \in \Ss$ is a product  of $\R_j \in \Rs$, the $X_i$ is a product of $\tl\R_j \in \tl\Rs$ with the same stoichiometric coefficient.
\end{definition}
We next provide a formal definition of the elementary modifications in Table I.
\begin{table}[t]
\centering
\includegraphics[width=0.585\columnwidth]{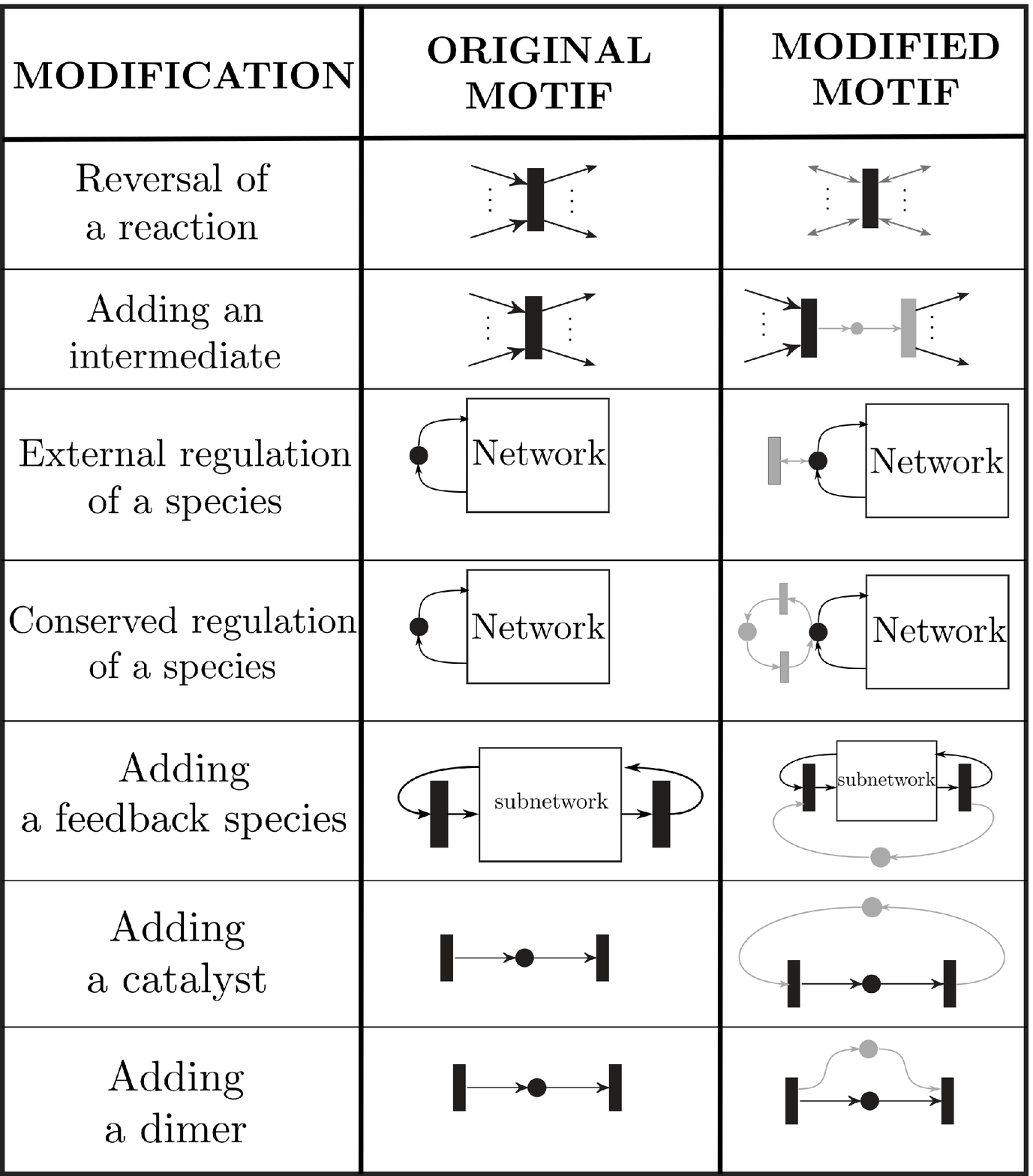}
\caption{A list of elementary graph modifications studied in this paper. Formal definitions are provided in Definition \ref{def.mod}.}
\end{table}
\begin{definition} \label{def.mod}
Let $\mathscr N=(\mathscr S,\mathscr R)$ be a given  BIN. We say that $\tl{\mathscr N}:=(\tl\Ss,\tl\Rs)$ is an \textit{elementary modification} of $\mathscr N$ if it satisfies one of the following statements:
\begin{enumerate}
	\item (\emph{Reversal of a reaction}) $\tl\Ss=  \Ss$, and $\exists \R_j \in \Rs$ such that $\tl\Rs=\Rs \cup \{\R_{-j}\}$. 
	\item (\emph{Adding an intermediate})  $\tl\Ss=\Ss \cup \{X^*\}$, and $\exists \R_j \in \Rs$ (written as $\R_j=\sum_{i}\alpha_{ij}X_i \to \sum_i \beta_{ij} X_i$) such that $\tl\Rs=(\Rs/\{\R_j\}) \cup \{\tl\R_{j},\tl \R^*\}$ where $\tl\R_j:= (\sum_{i}\alpha_{ij}X_i \to X^*)$, and $(\tl \R^*:=  X^* \to \sum_i \beta_{ij} X_i) $.
	\item (\textit{External Regulation}) $\tl\Ss=  \Ss$, and $\exists X_k \in \Ss$ such that $\tl\Rs=\Rs \cup \{ X_k \leftrightharpoons \emptyset \}$.
		\item (\textit{Conserved Regulation}) $\tl\Ss=  \Ss \cup \{X_{n+1}\}$, and $\exists X_k \in \Ss$ such that $\tl\Rs=\Rs \cup \{ X_k \leftrightharpoons X_{n+1} \}$.
	\item (\textit{Adding a feedback species}) $\tl\Ss=  \Ss \cup \{X^*\}$, and $\exists \R_j,\R_k \in \Rs$ such that $\tl\Rs=(\Rs/\{\R_j,\R_k\}) \cup \{\tl\R_{j},\tl \R_k\}$ where $\tl\R_j$ is an extension of $\R_j$ with $X^*$ as an extra product,  and $\tl\R_k$ is an extension of $\R_k$ with $X^*$ as an extra reactant.%
	\item (\emph{Adding a catalyst}) $\exists X_i \in \Ss$ such that 
	$\tl\Ss = \Ss \cup \{X_i^-\}$, $\vert\mathscr R\vert=\vert\tl\Rs\vert$, and every reaction  $\tl \R_j\in\tl\Rs$ is an extension of a corresponding reaction $\R_j \in \Rs$. Furthermore,  
	$X_i^-$ is a product   of a reaction $\R_j$ iff $X_i$ is a reactant of $\R_j$ with the same stoichiometry coefficient, and $X_i^-$ is a reactant of a reaction $\R_j$ iff $X_i$ is a product of $\R_j$ with the same stoichiometry coefficient. 
	
	\item (\emph{Adding a dimer})	 $\exists X_i \in \Ss$ such that				$\tl\Ss =  \Ss \cup \{X_i^+\}$, $\vert\mathscr R\vert=\vert\tl\Rs\vert$, and every reaction  $\tl \R_j\in\tl\Rs$ is an extension of a corresponding reaction $\R_j \in \Rs$. Furthermore, $\exists X_i \in \Ss$ such that  $X_i^+$ is a reactant of a reaction $\R_j$ iff $X_i^+$ is a reactant of $\R_j$ with the same stoichiometry coefficient, and $X_i^+$ is a product of a reaction $\R_j$ iff $X_i$ is a product of $\R_j$ with the same stoichiometry coefficient.
\end{enumerate} 
\end{definition}

Finally, a network $\tl \Ns$ is a modification of $\Ns$ if it is a result of several elementary modifications. More formally:
\begin{definition} A network $\tl\Ns$ is a \textit{modification} of $\Ns$ if there exists a finite sequence of networks $\Ns_0,\Ns_1,..,\Ns_q$, with $\Ns_0:=\Ns, \Ns_{q}:=\tl\Ns$, and for each $i\in{1,..,q}$, $\Ns_{i}$ is an elementary modification of $\Ns_{i-1}$. 
\end{definition}
In the subsequent sections, we provide results on modifications that preserve the stability of a given BIN. 

\begin{remark}\label{r.enzymatic} The standard enzymatic catalysis reaction is a combination of three elementary modifications which are  adding an intermediate, reversal, and then adding a catalyst. In other words, the reaction $S \to P$ is modified into $S \to C \to P$, then to $S \rightleftharpoons C \to P$, then to $S+E \leftrightharpoons C \to P+E$.
\end{remark}

\subsection{Linear networks with a Sum-of-Currents RLF}
It is easy to see that the first few modifications in Table I are stability preserving when applied to a linear BIN. This is stated below.
\begin{theorem}\label{th.linearmod1} Let $\Ns$ be a given \underline{linear} BIN, and let  $\tl\Ns$  be its modification generated by a finite sequence of elementary modifications that are limited to reversal of a reaction, adding an intermediate,  external regulation of a species, and conserved regulation of a species. Then, $V$ \eqref{e.lyap_linear} is an RLF for $\tl\Ns$.
\end{theorem}
\begin{proof} The resulting network $\tl\Ns$ after the application of the   elementary modifications mentioned in the statement of the theorem is linear. Hence, the statement follows by Theorem \ref{th1}.  
\end{proof}

The last two modifications in Table I are more interesting since they can modify a linear network into a nonlinear one. Nevertheless, we show that the resulting modified BIN continues to have an SoC RLF. The proof is provided in the appendix.

\begin{theorem} \label{mainth} Let $\Ns=(\Ss,\Rs)$ be a given \underline{linear} BIN, and let  $\tl\Ns=(\tl\Ss,\tl\Rs)$  be its modification generated by a finite sequence of elementary modifications that are limited to adding a catalyst and adding a dimer.    Then, $V=\sum_{i=1}^{\vert\mathscr S\vert}\vert\dot x_i\vert$ is an RLF for $\tl\Ns$.
\end{theorem}

Several modifications can be combined to yield enzymatic catalysis reactions (see Remark \ref{r.enzymatic}). Therefore, we can state the following corollary:
\begin{corollary}\label{cor} Let $\Ns$ be a given \underline{linear} BIN, and let  $\tl\Ns$  be its modification generated replacing linear reactions of the form $X_i \to X_j$, by nonlinear reactions of the form $X_i+E_{ij} \rightleftharpoons C_{ij} \to X_j + E_{ij}$. Let $\Ss_2$ be set of all the extra intermediates written as $C_{ij}$. 
Then, the function $V=\sum_{i=1}^{\vert\Ss\vert} \vert\dot x_i\vert + \sum_{C_{ij}\in\Ss_2} \vert\dot c_{ij}\vert$ is an RLF for $\tl\Ns$.
\end{corollary}
\begin{proof} The proof follows by using Theorem \ref{th.linearmod1} for adding an intermediate and then reversal, i.e., modifying $X_i\to X_j$ to $X_i \leftrightharpoons C \to X_j$. Then, Theorem \ref{mainth} to get the reaction $X_i +E_{ij}\leftrightharpoons C_{ij} \to X_j+E_{ij}$.
\end{proof}

\subsection{Networks with a Max-Min RLF}
Networks that have a Max-Min RLF admit a different set of stability-preserving modifications as we show next. Note that the original BIN does not need to be linear as is stated in the following result. 
\begin{theorem}\label{th.maxmin_mod} Let $\Ns$ be a BIN that admits a Max-Min RLF, and let  $\tl\Ns$  be its modification generated by a finite sequence of elementary modifications that are limited to adding an intermediate, adding a feedback species, adding a dimer, and adding a catalyst. Then, \eqref{maxmin} is an RLF for $\tl\Ns$.
\end{theorem}
\begin{proof} Using the characterization in Theorem \ref{th.maxmin}, any combination of the graph modifications mentioned in the statement of theorem do not create new independent vectors in the kernel of the stoichiometry matrix (i.e., it does not create new fluxes), and they do not make a single species a reactant in multiple reactions. Therefore, Theorem \ref{th.maxmin} applies to $\tl\Ns$.
\end{proof}

We study next the case of reversal. Since our formalism treats a reversible reaction as two reactions $\R_j,\R_{-j}$, then reversal of a reaction increases the number of fluxes, and hence violates the conditions required by Theorem \ref{th.maxmin}. Nevertheless, as shown in \cite{PWLRj}, the result can be extended. We state the result here in the language of graph modifications:
\begin{theorem}(\cite{PWLRj})\label{th.maxmin_reversal} Let $\Ns=(\Ss,\Rs)$ be a given network that satisfies the conditions of Theorem \ref{th.maxmin}. Let $\Rs_r \subset \Rs$ be defined as follows: $\R^* \in \Rs_r$ iff for each $X_i\in\Ss$ that is a product of $\R^*$, $X_i$ is not a product of another reaction. 
Then, let  $\tl\Ns$ be a modification of $\Ns$ generated by the reversal of the reactions in $ \Rs_r$. Then, \eqref{maxmin} is an RLF for $\tl\Ns=(\tl\Ss,\tl\Rs)$ where $\mathcal R(x)=\{R_{j}(x)-R_{-j}(x)\vert j=1,..,\vert\Rs\vert\}$, where $R_{-j}:\equiv 0$ if $\R_{-j}\not\in \tl\Rs $.
\end{theorem}

In addition, we can strengthen Corollary \ref{cor} to include modifications by \emph{processive} enzymatic cycles \cite{gunawardena07}:
\begin{corollary}\label{cor2} Let $\Ns$ be a given  BIN satisfying the conditions of Theorem \ref{th.maxmin}, and let  $\tl\Ns$  be its modification generated by replacing  reactions of the form $\sum_{i} \alpha_{ij} X_i \to \sum_{i} \beta_{ij} X_i$, by  reactions of the form $\sum_{i} \alpha_{ij} X_i+E^* \rightleftharpoons C_{0}^* \rightleftharpoons C_{1}^* \rightleftharpoons  .... \rightleftharpoons C_{m}^*  \to  \sum_{i} \beta_{ij} X_i + E^*$ for some positive integer $m$. Then, $\tl\Ns$ admits a Max-Min RLF.
\end{corollary}
\begin{proof}
	The statement can be proven by applying  enzymatic catalysis (as in Remark \ref{r.enzymatic}) to get $\sum_{i} \alpha_{ij} X_i+E^* \rightleftharpoons C_{0}^* \to \beta_{ij} X_i + E^*$, then by the addition of intermediates $C_1^*, .. C_m^*$ and then reversals to get the required reaction.
\end{proof}
\section{Global Stability and Robust Non-degeneracy}
\subsection{Global stability}
Since our RLFs are non-strict, we need to verify global stability. A popular way is via LaSalle's invariance principle. In our setting, a network $\Ns$ that  admits an RLF $V$ is said to satisfy the LaSalle's principle if the following implication always holds: If a bounded solution $\tl x(t)$ of \eqref{e.ode} satisfies $\tl x(t) \in \ker \dot V$ for all $t\ge 0$, then $\tl x(t) \in \ker V$ for all $t \ge 0$, i.e. $V(\tl x(t)))=0$. 

\begin{theorem}[\cite{PWLRj,MA_LEARN}] Let a network $\Ns$ be given. Assume that $\Ns$ admits an RLF and it satisfies the LaSalle's principle. Then, 
		\begin{enumerate} \item Each bounded trajectory converges to the set of steady states,
		\item if all the trajectories are bounded and there exists an isolated steady state relative to its stoichiometric class, then it is globally asymptotically stable.
	\end{enumerate} 
\end{theorem}
\subsubsection{Networks that admit an SoC RLF}
  In \cite{PWLRj}, an iterative algorithm has been proposed to check LaSalle's invariance principle. However, in the next result, we show that it always holds for networks that satisfy Theorems \ref{th.linearmod1} or \ref{mainth}.   
\begin{theorem}\label{th.global} Let $\mathscr N$ be a linear network or a modification of a linear network that satisfies the conditions of  Theorems \ref{th.linearmod1} or \ref{mainth}. Then it satisfies the LaSalle's principle.
\end{theorem}
 
For linear networks, the  statement  has been shown in \cite{maeda78}.  It remains to prove that generalization of the the result to any nonlinear network that is a modification of a linear network. The proof is included  in the Appendix. 

\subsubsection{Networks that admit a Max-Min RLF}
Verification of LaSalle's invariance principle  for networks that admit Max-Min RLFs has been provided in \cite{PWLRj} via a simple graphical condition. In order to introduce the next result, we need a definition. Consider a network $\Ns=(\Ss,\Rs)$, then a reaction $\R_k \in \Rs$ is said to be an ancestor of $\R_j$ if there is a direct path from $\R_k$ to $\R_j$ on the Petri-net corresponding to $\Ns$.  The result is stated in the following theorem:
\begin{theorem}[\cite{PWLRj}] \label{th.maxmin_ls} Let $\mathscr N$ be a  network that satisfies the conditions of Theorem \ref{th.maxmin}, Theorem \ref{th.maxmin_mod}, or Theorem \ref{th.maxmin_reversal}.  Then, $\Ns$  satisfies the LaSalle's principle: if $\Ns$ is conservative, or if  every pair of reactions share an ancestor.
	\end{theorem}

\subsection{Robust Non-degeneracy}
\subsubsection{Definitions and review}
In the previous subsection, we have shown that the trajectories converge to the set of steady states. However, existence of a steady state in a stoichiometric class does not automatically imply that it is isolated. Therefore, we study here the robust non-degeneracy of the Jacobian of \eqref{e.ode} which can be written as $\Lambda:=\Gamma  \partial R/\partial x =\Gamma  V$, where $V\in \mathcal K_\Ns $.    However, as mentioned in §2, the presence of a conservation law means that the positive orthant is a foliage of forward invariant stoichiometric classes. Therefore, the relevant entity for analysis is the \emph{reduced Jacobian} $\Lambda_r$ which can be defined as follows. For a given $\Gamma \in\mathbb R^{n\times \nu}, V \in\mathbb R^{\nu\times n}$, denote $r:=\mbox{rank}(\Gamma)$.  Let $\{d_1,..,d_{n-r}\}$ be linearly independent left null vectors of $\Gamma$. In order to get a basis of $\mathbb R^n$, we add  vectors to get the basis: $\{d_1,..d_{n-r},d_{n-r+1},..,d_{n}\}$, and get the transformation matrix:
\[ T = [ d_1,...,d_n]^T.\]
The Jacobian $\Lambda=\Gamma V$ in the new coordinates can be written as follows:
\begin{equation}
	\label{r.J} T \Gamma V T^{-1} = \begin{bmatrix} \Lambda_r & \Lambda_2 \\ 0 & 0 \end{bmatrix} 
\end{equation}
The matrix $\Lambda_r \in \mathbb R^{r \times r}$ is the \emph{reduced Jacobian}, and it is the Jacobian for the dynamics restricted to the stoichiometric class.

 We are interested in its non-singularity for any admissible kinetics.  Hence, we provide the following definition:
\begin{definition}
	A network $(\Ss,\Rs)$  is said to be \emph{robustly non-degenerate} iff the reduced Jacobian $\Lambda_r$ defined in \eqref{r.J} is non-singular for all matrices $V \in \mathcal K_\Ns $.
\end{definition}

Although \eqref{r.J} is written with a specific transformation matrix $T$, it is obvious to see that the non-singularity of the reduced Jacobian  is independent of the specific choice of the matrix $T$.

In order to study the reduced Jacobian, we will use the concept of the \emph{essential determinant} $\mbox{det}_{{ess}}(\Lambda)$  which is defined as the sum of all $r\times r$ principal minors of $\Lambda$. The characterization can be stated as follows:
\begin{lemma}(\cite{banaji16}) \label{l.essdet} Let $\Gamma \in \mathbb R^{n\times \nu}$ and $V \in \mathbb R^{\nu \times n}$ be given. Let $r:=\mbox{rank}(\Gamma)$. The reduced Jacobian $\Lambda_r$ defined {in} \eqref{r.J} is non-singular iff $\det_{ess}(\Lambda)=\det_{ess}(\Gamma V)\ne 0$.
	\end{lemma}

Hence, instead of explicitly computing the reduced Jacobian, our strategy will be to verify that the sum of $r\times r$ {principal} minors of the full Jacobian is nonzero for any   $V \in \mathcal K_\Ns $. 
Our task is eased by the special properties of  networks admitting a PWL RLF. In \cite{PWLRj}, we have proved  that every principal minor of the Jacobian is, in fact, non-negative. To state it more formally, we have the following definition:
\begin{definition} A network $\Ns$  is said to be \emph{robustly $P_0$} if
	the Jacobian $-\Gamma V  $ is $P_0$ for all $V\in\mathcal K_\Ns$, i.e., all its principal minors are non-negative.
\end{definition}

Hence, the result can be stated as follows.
\begin{lemma}(\cite{PWLRj}) \label{lem.p0} Let $\Ns$ be a given network. If it admits a PWL RLF, then it is robustly $P_0$.
	\end{lemma}

Therefore, using the last two lemmas, we immediately get the following corollary:
\begin{corollary}\label{cor.nondeg} Let $\Ns$ be a given network that is robustly $P_0$. Then, $\Ns$ is robustly non-degenerate iff for every $V\in \mathcal K_N$, there exists a positive $r\times r$ principal minor, where $r=\rank(\Gamma)$.
\end{corollary}

\subsubsection{Computational testing of robust non-degeneracy}
It is possible to computationally check robust non-degeneracy by testing the Jacobian at a finite number of points \cite{blanchini14,MA_LEARN,dissertation}. In fact, we have shown that one point is sufficient:
	\begin{theorem}[\cite{dissertation,MA_LEARN}] \label{th15} Let $\Ns$ be a network that admits a PWL RLF, and let $\G\in\mathbb R^{n\times \nu}$ be the stoichiometry matrix with rank $r${.} If $\exists V^* \in \mathcal K_{\Ns} $ such that $-\G V^*$ has a positive essential determinant, then $-\G V^*$ has a positive essential determinant for all $V\in \mathcal K_\Ns $, i.e,  $\Ns$ is robustly non-degenerate.  
	\end{theorem}

 In the next subsection, we provide our main result in this subsection, which is a graphical method to verify robust non-degeneracy.
 
 \begin{remark} Although  Theorem \ref{th15} is stated in \cite{MA_LEARN,dissertation} for networks admitting a PWL RLF, the proof holds for any robustly $P_0$ network. \end{remark}
\subsubsection{Main Result}
Instead of directly verifying the non-degeneracy of a large network, we study it graphically. In other words, we consider the network as a modification of a simpler network. Therefore, we state our result which is proved in the appendix.
\begin{theorem}\label{th.linearmod_J} Let $\Ns$ be a given   BIN which is robustly $P_0$. Assume that $\Ns$ is robustly non-degenerate, and let  $\tl\Ns$  be its modification generated by a finite sequence of elementary modifications that are limited to reversal of a reaction, adding an intermediate,  external regulation of a species, conserved regulation of a species, adding a catalyst, and adding a dimer. Then, if $\tl\Ns$ is robustly $P_0$, it follows that $\tl\Ns$ is robustly non-degenerate.
\end{theorem}
\begin{remark}
Any network that admits a PWL RLF is robustly  $P_0$ using Lemma \ref{lem.p0}. Hence, Theorem \ref{th.linearmod_J} can be coupled with Theorems \ref{th.linearmod1},\ref{mainth},\ref{th.maxmin_mod},\ref{th.maxmin_reversal} to show robust non-degeneracy of the modified networks.
\end{remark}

\subsection{Review of the consequences of robust non-degeneracy}
Robust non-degeneracy of the Jacobian gives us a quick way to verify several key properties of BINs.  For completeness, we review them here.

\paragraph{Uniqueness of steady states} Lemma \ref{lem.p0} implies that any network that admits a PWL RLF has a $P_0$ Jacobian, {which excludes multiple non-degenerate steady states in the same stoichiometric class \cite{banaji07,banaji16}}. Hence, we get the following:
\begin{theorem}[\cite{dissertation,MA_LEARN}] Consider a network $\Ns$ that admits a PWL RLF and is robustly non-degenerate. Then every positive steady state is unique relative to its stoichiometric class.
	\end{theorem}

\paragraph{Exponential stability} The following result follows from the properties of PWL Lyapunov functions:
\begin{theorem}[\cite{blanchini14,dissertation,MA_LEARN}] Let $\Ns$ be a network that admits a PWL RLF and robustly non-degenerate, then every positive steady state  is exponentially asymptotically stable.
	 \end{theorem}

 \paragraph{Global stability} Using previous results, it can be readily seen that robust non-degeneracy coupled with the LaSalle's principle {implies} global stability. However, it has been shown \cite{blanchini17} that this can be strengthened to the following:
 \begin{theorem}[\cite{blanchini17}] Suppose that a network $\Ns$ admits a PWL RLF and is robustly non-degenerate, then every positive steady state  is globally asymptotically stable relative to its stoichiometric class.
 \end{theorem}

\begin{remark}
	The statements in this subsection assume the existence of a \emph{positive} steady state. One way to exclude the existence of steady states on the boundary is via verifying \emph{persistence}. In other words, we need to guarantee that all the trajectories that start from the positive orthant do not asymptotically approach its boundary. Graphical conditions for persistence have already been developed in \cite{angeli07p} and they are easily applicable as we will see in the next section.
\end{remark}

{
\section{Persistence}
\subsection{Definitions and review of previous results}
For systems that evolve on the positive orthant, persistence simply means \emph{non-extinction} \cite{gard80,waltman91}. In other words, if a trajectory starts in the interior of the positive orthant, then it will not approach the boundary asymptotically. More precisely, a trajectory   $\varphi(t;x_\circ)$ of \eqref{e.ode} is said to be \textit{persistent} if it satisfies $\liminf_{t\to\infty}\varphi(t;x_\circ)\gg0 $ whenever $x_\circ\gg 0$. A BIN network $\Ns$ is said to be \emph{robustly persistent} if the previous statement holds for all bounded trajectories and for all admissible kinetics. A graphical notion of robust persistence  for BINs has been introduced in \cite{angeli07p,Angeli07} using the concept of \emph{siphon} which we define next. \begin{definition} Let $\Ns=(\Ss,\Rs)$ be a given network. Then,   a non-empty set $P \subset \Ss$ is called a \emph{siphon} iff each input reaction associated to a species in $P$ is also an output reaction associated to a (possibly-different) species in $P$. 
		 A siphon is said to be \emph{trivial} if it contains the support of a conservation law, and it is said to be \emph{critical} otherwise.
	 
\end{definition}
The main result is as follows:
\begin{theorem}(\cite{Angeli07,angeli07p,angeli11})
	Let $\Ns$ be   a network that lacks critical siphons. Then, $\Ns$ is robustly persistent.
\end{theorem}
This motivates the following definition:
\begin{definition}
A network $\Ns$ that lacks critical siphon is said to be \emph{graphically persistent}.
\end{definition}

\subsection{Main Result}
We show here that graphical persistence is conserved under many types of modifications. We start with a general result whose proof is provided in the appendix:
\begin{theorem}\label{th_cs1} Let $\mathscr N$ be a given  BIN which is graphically persistent, and let $\tl\Ns$ be its modification generated by a finite sequence of elementary modifications that are limited to reversal of a reaction, external regulation of a species, conserved regulation of a species, adding an intermediate, and adding a dimer. Then, $\tl\Ns$ is also graphically persistent.
	\end{theorem}
We next show that the last theorem can be expanded for the classes of networks studied in this paper:
\begin{theorem}\label{th_cs2} Let $\mathscr N$ be a given  BIN which is graphically persistent.  
	\begin{enumerate}
		\item Let $\Ns$ be linear, and let $\tl\Ns$ be its modification using any of the modifications listed in Theorems \ref{th.linearmod1} or \ref{mainth}, then $\Ns$ is graphically persistent.
		\item Assume that $\Ns$ is conservative and that it admits a Max-Min RLF, then it is graphically persistent.
	\end{enumerate} 
\end{theorem}
Part 1 of Theorem \ref{th_cs2} follows from Theorem \ref{th_cs1} except for the case of adding a catalyst which is proved in the Appendix. Part 2 of Theorem \ref{th_cs2} is proved in \cite[Theorem 13]{PWLRj}.
}

\section{Applications}
\subsection{Post-translational Modification (PTM) cycles}
The PTM cycle model  is standard in systems biology  \cite{goldbeter81}.
The long-term dynamics of the PTM cycle have been a subject of extensive study using several methods. This includes monotonicity \cite{angeli08,angeli10}, and RLFs \cite{PWLRj,MA_LEARN,blanchini14}. In this paper, we show that the stability properties of the PTM cycle can be interpreted graphically in terms of the basic reversible reaction:
\begin{equation}\label{reversible}S \rightleftharpoons P,\end{equation}
where $S$ denotes the substrate, and $P$ denotes the product. This simple motif admits both an SoC RLF and a Max-Min RLF. Furthermore, it is conservative, robustly non-degenerate, satisfies the LaSalle's condition. In addition, it lacks critical siphons, hence it is persistent \cite{angeli07p}. Therefore, it satisfies the following statement which we call $(\star)$: \textit{Each proper stoichiometric class contains a unique globally exponentially stable positive steady state}. We show next that these properties are inherited by the modifications of the simple reversible reaction above.

\subsubsection{The single PTM} We consider the single PTM cycle:
\begin{equation} S + E \leftrightharpoons C \lra P+E,   P + F \leftrightharpoons D \lra S+E.  \end{equation}
As noted in Remark \ref{r.enzymatic}, the reaction $S \to P$ can be modified into an enzymatic catalysis reaction.   Using Corollary \ref{cor}  we  get that the PTM cycle above admits an SoC RLF. Furthermore, using Corollary \ref{cor2} we get that it also admits a Max-Min RLF. {Theorem \ref{th_cs2} implies that it lacks critical siphons.} Hence, using the results in §5  it satisfies the statement $(\star)$.

\subsubsection{The PTM star}\label{s.ptmstar}
\begin{figure}
\centering
\includegraphics[width=0.983586257\columnwidth]{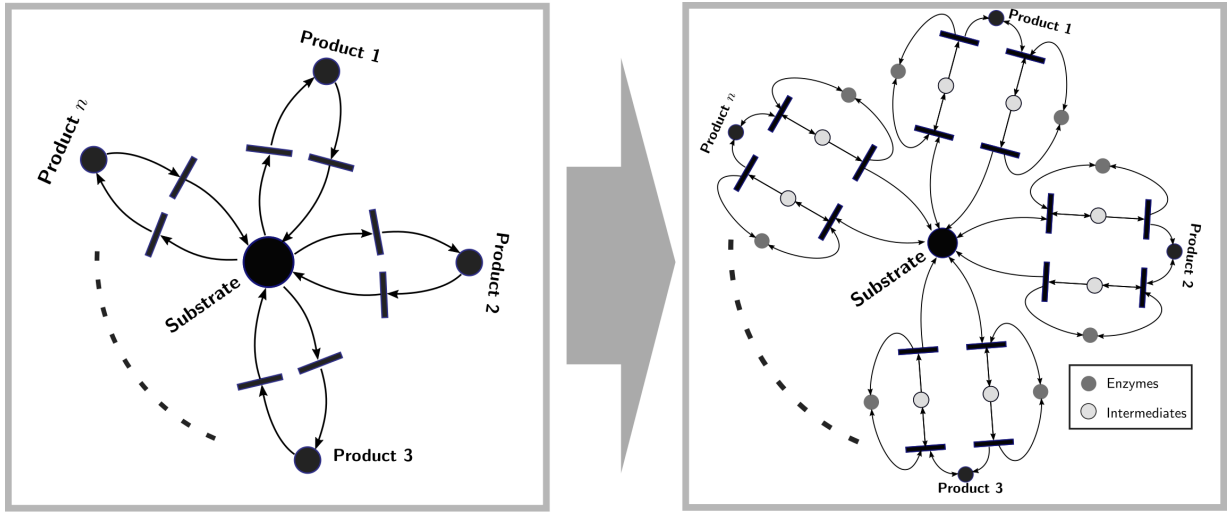}
\caption{\textbf{Stability of the linear star implies stability of the PTM star}. Using Corollary \ref{cor}, existence of an RLF for the linear star implies the existence of an RLF for the PTM star depicted in Fig. \ref{f.star}.}
\label{fig:ptmstarlinear}
 
\end{figure}

We can consider other modifications to \eqref{reversible}. By adding a finite number of conserving regulations on $S$, we get the following network which we call the linear star (depicted in Fig. \ref{fig:ptmstarlinear}):
\begin{equation}\label{star_linear} S \rightleftharpoons P_1,  S \rightleftharpoons P_2,..., S \rightleftharpoons P_n. \end{equation}
Then, using Corollary \ref{cor}, we get that the PTM star \eqref{starA}-\eqref{starB} (depicted in Fig. \ref{f.star}) admits an SoC RLF since it is formed by enzymatic catalysis modifications. {In addition, Theorem \ref{th_cs2} implies that it lacks critical siphons.}. Hence, it satisfies that statement $(\star)$. Furthermore, to answer the question posed in the introduction. We can add the external regulation $\emptyset \leftrightharpoons S$ to \eqref{star_linear}, and then apply enzymatic catalysis to all other reactions to certify the existence of an SoC RLF. Since the network is no longer conservative, it satisfies the following statement: \textit{If a proper stoichiometric class contains a steady state, then it is a unique globally exponentially stable positive steady state.}

\subsubsection{The Processive Multi-PTM  cycle.} Modifying \eqref{reversible} by adding intermediates gives the following network which we call the linear cycle (depicted in Fig. \ref{f.chain_cycle}-a):
\[ S_0 \to S_1 \to ... \to S_n \to S_0,\]
where $S_0:=S, S_n:=P$. Theorem \ref{th.maxmin} guarantees that the modified network has a Max-Min RLF.  Corollary \ref{cor2} implies that the following network admits a Max-Min RLF:
\begin{align*}  
S_{i-1}+E_i & \rightleftharpoons C_{i1} \rightleftharpoons C_{i2} \rightleftharpoons .... \rightleftharpoons C_{im} \longrightarrow S_i + E_i,  \\  
S_{n}+E_n & \rightleftharpoons C_{n1} \rightleftharpoons C_{n2} \rightleftharpoons .... \rightleftharpoons C_{nm} \longrightarrow S_1 + E_n,
\end{align*}
$i=1,..,n-1$. The above network has been called the ``all-encompassing'' processive cycle, and its stability  has been studied in \cite{shiu17} using {monotone system techniques}. Using our method, we show that the existence of an RLF follows by modifying the linear cycle (Fig. \ref{f.chain_cycle}-a) using processive enzymatic reactions to get the network depicted in Fig. \ref{f.chain_cycle}-b. {In addition, Theorem \ref{th_cs2} implies that it lacks critical siphons.} Therefore, using the results in §5, it satisfies the statement $(\star)$.

\subsubsection{The PTM chain}  
\begin{figure}
\includegraphics[width=\columnwidth]{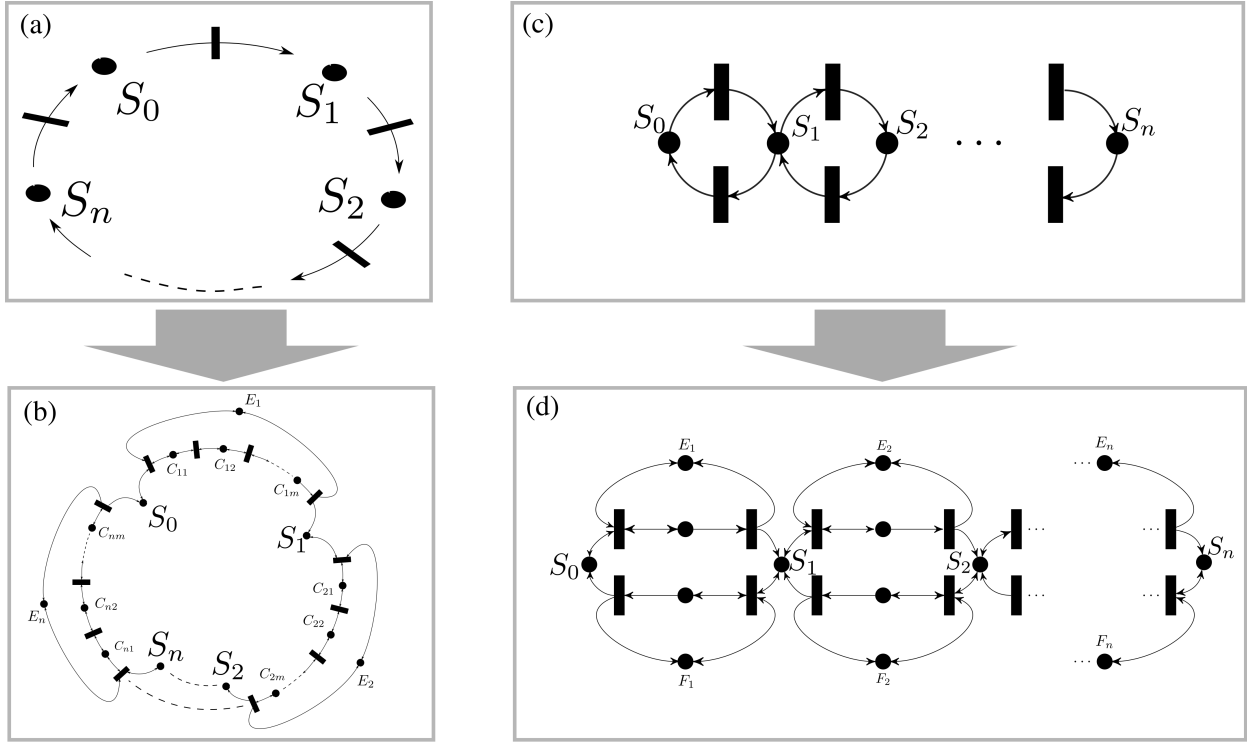}
\caption{\textbf{Constructing an RLF for a nonlinear network from a linear one.} (a) The linear cycle. (b) The processive multi-PTM cycle. The existence of an RLF follows from the existence of one for the linear cycle using Corollary \ref{cor2}. (c) The linear chain. (d) The PTM chain. The existence of an RLF follows from the existence of one for the linear chain using Corollary \ref{cor}.  }
\label{f.chain_cycle}
\end{figure}

Consider now modifying \eqref{reversible} by  a finite number of intermediates and reversals, we get the following network:
\begin{equation}\label{chain_linear} S_0 \rightleftharpoons S_1 \rightleftharpoons S_2 ... \rightleftharpoons S_n, \end{equation}
where $S_0:=S, S_n:=P$.
Corollary \ref{cor} implies that the following PTM chain admits an SoC RLF:
\begin{align}  
S_{i-1}+E_i & \rightleftharpoons C_{i} \longrightarrow S_i + E_i, \\  
S_{i-1}+F_i & \rightleftharpoons D_i \longrightarrow S_i + F_i, i=1,..,n.
\end{align}
The existence of an SoC RLF of the PTM chain can be shown computationally for each given $n$ by linear programming \cite{MA_LEARN}. Nevertheless,
Fig. \ref{f.chain_cycle}-c,d shows that the existence of an SoC RLF \emph{for each $n$} follows from modifying a linear  chain via enzymatic catalysis reactions.   {In addition, Theorem \ref{th_cs2} implies that it lacks critical siphons.} Therefore, using the results in §5, it satisfies the statement $(\star)$.

\begin{figure}
\centering
\includegraphics[width=0.95\columnwidth]{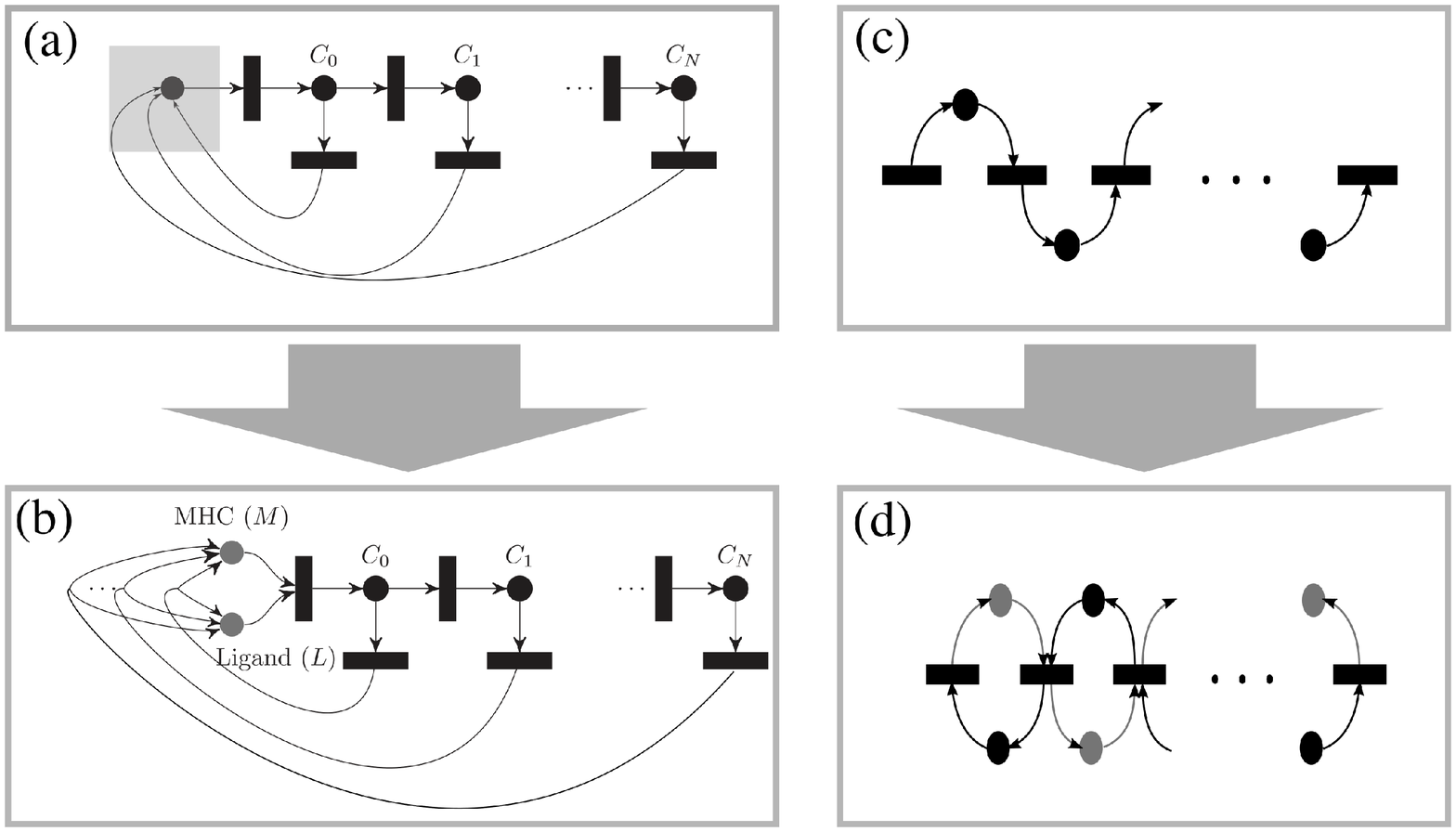}
\caption{\textbf{Additional examples for graphical RLF construction.} (a) A linear BIN. (b) The McKeithan network.  The existence of an RLF follows from the existence of one for the linear BIN in panel (a) using Corollary \ref{cor}. (c) A one-directional linear chain. (d) The RFM. The existence of an RLF follows from the existence of one for the unidirectional linear chain using Corollary \ref{cor}.  }
\label{f.tcell_rfm}
 
\end{figure}

\subsection{$T$-cell kinetic proofreading}
McKeithan \cite{mckeithan95} proposed a nonlinear BIN to explain $T$-cell's ability to distinguish between different types of ligands. It is given as follows:
\begin{align} \label{mckeithen_nl}
R+L &\rightleftharpoons C_0 \to C_1 \to ... \to C_n \\ \nonumber
C_1 &\to R+L, C_2 \to R+L,..., C_n \to R+L.
\end{align}

Sontag \cite{sontag01} has studied the stability of the network using the theory of complex balance, while we have studied the network using computational RLF construction \cite{MA_LEARN}. Here, we show that a stability certificate can be constructed by considering the network as a modification of a linear network. 
By noting that the species $L$ is a dimer in the language of Table I, we can see that \eqref{mckeithan} is a modification of the following network by the addition of a dimer:
\begin{align}\label{mckeithan}
RL &\rightleftharpoons C_0 \to C_1 \to ... \to C_n \\ \nonumber
C_1 &\to RL, C_2 \to RL,..., C_n \to RL.
\end{align}
 Hence, existence of an SoC RLF for \eqref{mckeithan} follows from Theorem \ref{mainth}. Fig. \ref{f.tcell_rfm}-a shows the linear network, while Fig. \ref{f.tcell_rfm}-b shows the corresponding modified nonlinear network.

The set of steady states is globally stable by  Theorem \ref{th.global}. We can also show robust non-degeneracy graphically as follows. We consider first a linear cycle $RL \to C_0 \to .. \to C_n \to RL$ which is robustly non-degenerate since it is a modification of $RL \rightleftharpoons C_n$. Then, adding   reactions of the form $C_i \to RL$ won't increase the rank of the stoichiometry matrix, hence the network in \eqref{mckeithan} is robustly non-degenerate using the same argument used in the proof of item 1 in Theorem \ref{th.linearmod_J}. Finally, \eqref{mckeithen_nl} is a modification of \eqref{mckeithan} by the addition of a dimer. Hence, robust non-degeneracy of \eqref{mckeithen_nl} follows from Theorem \ref{th.linearmod_J}.   {In addition, Theorem \ref{th_cs2} implies that it lacks critical siphons.} Therefore, using the results in §5, it satisfies the statement $(\star)$ for any $N$.

\subsection{The Ribosome Flow Model}
The Ribosome Flow Model (RFM) is a nonlinear system model of the process of translation initiation and elongation where it describes Ribosome binding to codons on an mRNA that is being translated \cite{rfm11}. It has been shown \cite{MA_LEARN} that the corresponding ODE can be written as a BIN with species $X_i, Y_i$ where $X_i$ is \emph{occupancy} of the $i$th codon, while $Y_i$ is the \emph{vacancy} of the $i$th codon.  Hence, we get the following BIN (depicted in Fig. \ref{f.tcell_rfm}-d):
\begin{align*} Y_1 &\to X_1, X_n \to Y_n, \\ X_i+Y_{i+1} &\to X_{i+1}+Y_{i}, \, i=1,..,n-1. 
\end{align*}
The stability of the above network has been studied via monotonicity methods \cite{rfm12}. For a given $n$, the existence of an SoC RLF can be verified via linear programming \cite{MA_LEARN}. Nevertheless, Fig. \ref{f.tcell_rfm}-c,d shows that an SoC RLF can be constructed by merely noticing that the RFM is a modification generated by adding catalysts to the following unidirectional linear chain network: (depicted in Fig. \ref{f.tcell_rfm}-c)
\[ \emptyset \to X_1 \to X_2 \to ... \to X_n \to \emptyset. \]

\begin{figure}
\centering
\includegraphics[width=0.635\columnwidth]{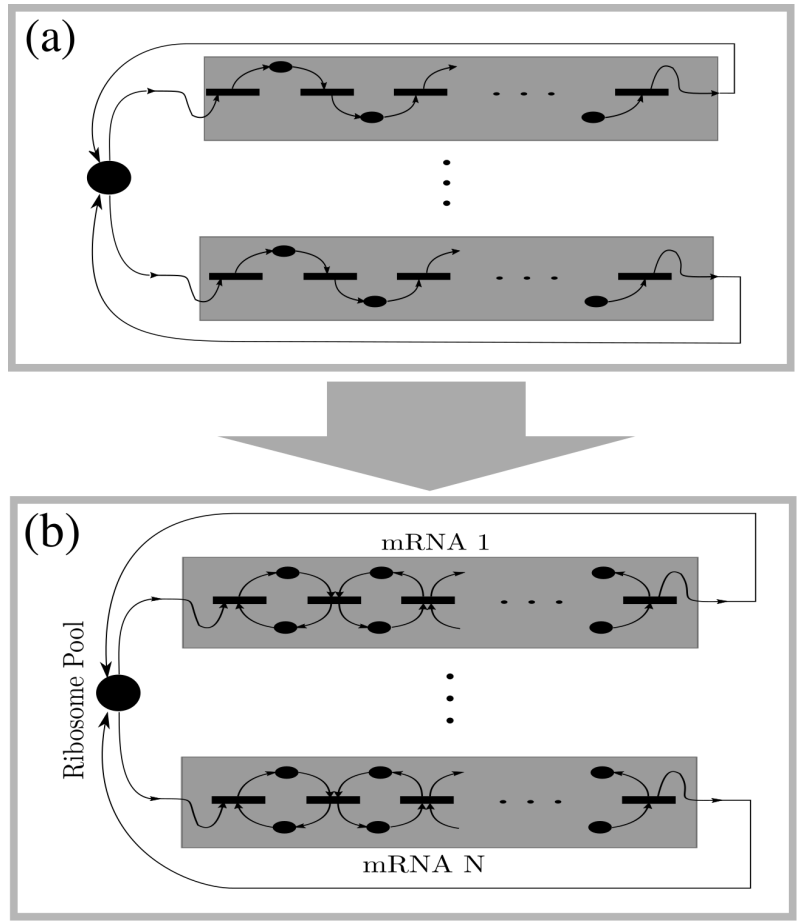} 
\caption{\textbf{Graphical construction of an RLF for the RFM with a pool \cite{rfm16}}. (a) The linear network. (b) The corresponding modified nonlinear network. Stability follows via Corollary \ref{cor}. \label{f.rfmpool}}
 
\end{figure}

The same graphical technique can be applied to RFMs interconnected via a pool \cite{rfm16} (as Figure \ref{f.rfmpool} shows), or via multiple pools \cite{LCSS_jared}.  {In addition, Theorem \ref{th_cs2} implies that they lack critical siphons.} Therefore, using the results in §5, and by noting that they lack critical siphons, all the aforementioned RFM variants satisfy the statement $(\star)$.%
\section{Conclusion}
In this work, we have proposed a graphical method to certify the existence of an RLF for a given network by reducing it via a certain set of admissible modifications to a network that is known to admit an RLF. Furthermore, our method can directly show that the stability of a given network is preserved under certain graph modifications. In addition, we have shown that  properties of the original network such as global stability, robust non-degeneracy, and {graphical persistence} are invariant under such modifications.  {Using our methods, complex nonlinear networks of arbitrary size and arbitrary number of nonlinear reactions can be reduced into tractable    networks.}  %

\section*{Appendix: Proofs}
\subsubsection*{Proof of Theorem \ref{th1}}

	The function $V(x)=\tV(R(x))$ is piecewise linear in terms of the rates, therefore there exists a positive integer $m$ such that the space $\mathbb R_{\ge0}^\nu$ can be partitioned into non-empty-interior regions $\{\cW_k\}_{k=1}^m \subset \mathbb R_{\ge0}^\nu$ for which $\tV$ is linear on each of them and each region corresponds to a specific sign pattern for $\dot x$. The geometry of such partition is discussed more thoroughly in \cite{PWLRj}.
	
	Fix $k$. There exists $c_{ij}^{(k)}, c_i^{(k)}, \delta^{(k)}$ such that: \begin{align}\label{p1} V(x)&=\sum_{\substack{i,j\\ i\ne j }} c_{ij}^{(k)} R_{ij}(x_i) + \sum_{i} c_{i}^{(k)} R_{i}(x_i) + c_0^{(k)} \\ \nonumber & =:c^{(k)^T}R(x), \ R(x)  \in \cW_k. 
	\end{align}
	Since $V$ is defined as the $\ell_1$ norm of  $\dot x$, then the sign of $\dot x$ is constant and non-zero on $\cW_k^\circ$. Therefore, we denote $\sigma_i:=\sgn(\dot x_i) \in \{\pm 1\}$ on $\cW_k^\circ$, where the superscript ``$\circ$'' denotes the interior of a set.
	
	We claim that each term in the expression \eqref{p1} has a nonpositive Lie derivative on $\cW_k^\circ$. In order to show that, we first examine terms of the form $c_{ij}^{(k)} R_{ij}(x_i)$ where $c_{ij}^{(k)} \ne 0$ for some $i,j$. We will show that $ c_{ij}^{(k)} \dot R_{ij}(x_i)  \le 0$ for $R(x) \in \cW_k^\circ$.   As evident by examining \eqref{e.lyap_linear}, the reaction rate $R_{ij}$ appears only in  $\dot x_i$ with coefficient $-1$ and in $\dot x_j$ with coefficient $+1$. W.l.o.g, assume that $c_{ij}^{(k)}>0$. There are four possible combinations $\sigma_i^{(k)},\sigma_j^{(k)}>0$, $\sigma_i^{(k)},\sigma_j^{(k)}<0$, $\sigma_i^{(k)}>0,\sigma_j^{(k)}<0$, and $\sigma_i^{(k)}<0,\sigma_j^{(k)}>0$. The first two give $c_{ij}^{(k)}=0$ and the third gives $c_{ij}^{(k)}=-2<0$. Hence, we conclude that  $\sigma_i^{(k)}<0,\sigma_j^{(k)}>0$.  %
	Therefore, $ \sgn(c_{ij}^{(k)} \dot R_{ij}(x_i))\! =\! \sgn(c_{ij}^{(k)} (\partial R_{ij}(x_i)/\partial x_i)   \dot x_i) \mathop{=}\limits^{ } \sgn(c_{ij}^{(k)} \sigma_i^{(k)})\le 0$ for $R(x) \in \cW_k^\circ$, where the last equality follows by the monotonicity of $R_{ij}$.
	
	Next, we examine $c_i^{(k)} R_i(x_i)$ for some $i$ where $c_i^{(k)} \ne 0$. W.l.o.g, assume that $c_i^{(k)}>0$. Since $R_i$ appears only in $\dot x_i$ with coefficient $-1$, then $\sigma_i^{(k)}<0$. Therefore, $ \sgn(c_{i} \dot R_{i}(x_i)) = \sgn(c_{i}^{(k)} (\partial R_{i}(x_i)/\partial x_i)   \dot x_i) \mathop{=}\limits^{ } \sgn(c_{i}^{(k)} \sigma_i)\le 0$ for $R(x) \in \cW_k^\circ$.   
	
	Since $k, i,j$ have been chosen arbitrarily, we conclude that $\dot V(x) \le 0$ whenever $R(x) \in \cW_k^\circ$ for some $k$. It remains to show that $\dot V(x) \le 0$ when $R(x) \in \partial \cW_k$ for some $k$ where ``$\partial$'' denotes the boundary of  a set. To that end, similar to \cite{PWLRj}[Proof of Theorem 2], the Dini's derivative can be written as $\dot V(x)=\max_{k \in K_{x(t)}} {c^{(k)}}^T \dot R(x) \le 0$ where $K_{x(t)}=\{k \vert R(x) \in \cW_k\}$. \\ \strut
\hfill $\square$

\subsubsection*{Proof of Theorem \ref{mainth}}
 
	Let $\Gamma$ be the stoichiometry matrix for $(\Ss,\Rs)$. Since the modifications are limited to adding a catalysis or adding a dimer, then every reaction in $\tl\Rs$ is an extension of a corresponding reaction in $\Rs$.  Hence, we can write $\tilde\Gamma=[\Gamma^T,\Gamma_2^T]^T$ as the stoichiometry matrix for $(\tilde\Ss,\tilde\Rs)$. Let $\dot x=\Gamma R(x), \dot{\tilde x}=\tilde\Gamma \tilde R(\tilde x)$ be the corresponding ODEs. Hence, we can write $\tilde x=[x^T,x_2^T]^T$, where $x_2$ corresponds to the concentrations of the species in $\tilde\Ss/\Ss$. 
	
	Note that all the  species in $\tilde\Ss/\Ss$ are either catalysts or dimers.  We include an additional assumption to simplify the notation:  For each species $X_i \in \Ss$, we assume that there exists at most one corresponding catalyst species in $\tilde\Ss/\Ss$, and it is denoted by $X_i^-$. Similarly, we assume that there exists at most one corresponding dimer species, and the corresponding species is denoted as $X_i^+$. The corresponding concentrations are $x_i, x_i^-, x_i^+$.  
	The proof can be generalized easily without the last assumption . 
	
	Our construction implies that $\dot x_i=\dot x_i^+ = -\dot x_i^-$. Hence, $V(x)=0$ iff $\dot x=0$. Therefore, $V$ is positive-definite. We next show that it is non-increasing. 
	
	Similar to the proof of Theorem \ref{th1}, we consider a region $\cW_k$ for which $V$ is linear and has a fixed sign pattern for $\dot x$. Fix $k$, There exists $c_{ij}^{(k)}, c_i^{(k)}, \delta^{(k)}$ such that: 
	\begin{align}\label{p2} V(x)&=\sum_{\substack{i,j\\ i\ne j }} c_{ij}^{(k)} R_{ij}(x_i,x_i^+,x_j^-) + \sum_{i} c_{i}^{(k)} R_{i}(x_i,x_i^+) + c_0^{(k)} \\ \nonumber & =:c^{(k)^T}R(x), \ R(x)  \in \cW_k. 
	\end{align}
	
	We claim that each term in the expression \eqref{p1} has a nonpositive Lie derivative on $\cW_k^\circ$. In order to show that, we first examine  $c_{ij}^{(k)} R_{ij}(x_i)$ where $c_{ij}^{(k)} \ne 0$ for some $i,j$. We will show that $ c_{ij}^{(k)} \dot R_{ij}(x_i)  \le 0$ for $R(x) \in \cW_k^\circ$.   Since  the candidate RLF sums only the species in $\Ss$, the reaction rate $R_{ij}$ appears only in  $\dot x_i$ with coefficient $-1$ and in $\dot x_j$ with coefficient $+1$. W.l.o.g, assume $c_{ij}^{(k)}>0$. Similar to the proof of Theorem 1, we get that $\sigma_i^{(k)}<0,\sigma_j^{(k)}>0$.
	Since $\sigma_i^{+(k)} = \sigma_i^{(k)}$, and  $\sigma_j^{-(k)} = -\sigma_j^{(k)}$ we  get \begin{align*}  \sgn(c_{ij}^{(k)} \dot R_{ij}(x_i,x_i^+,x_j^-))    &=  \sgn   \left (c_{ij}^{(k)} \left (\frac {\partial R_{ij} }{\partial x_i}    \dot x_i + \frac {\partial R_{ij} }{\partial x_i^+}    \dot x_i^+  + \frac {\partial R_{ij} }{\partial x_j^-}    \dot x_j^-   \right ) \right ) \\ & \mathop{=}\limits^{ } \sgn(c_{ij}^{(k)} (\sigma_i^{(k)}+\sigma_i^{+(k)}+\sigma_j^{-(k)} ))\le 0 \end{align*} for $R(x) \in \cW_k^\circ$, where the last equality follows by the monotonicity of $R_{ij}$.
	
	Since $k, i,j$ have been chosen arbitrarily, we can use the same arguments used in the proof of Theorem 1 to conclude that $\dot V(x)\le 0$ for all $x$.\\ \strut
\hfill $\square$

\subsubsection*{Proof of Theorem \ref{th.global}}
Let $\tl\Ns=(\tl\Ss,\tl\Rs)$ be a modification of a linear network $\Ns=(\Ss,\Rs)$ by adding a catalysts or   dimers. We use the standard LaSalle's principle \cite{yoshizawa}. Let $x(t)$ be a trajectory of \eqref{e.ode} that is contained in $\ker \dot V$, i.e., $\dot V(x(t))\equiv 0$, and $V(x(t))\equiv \bar V_1 \ge 0$.  In order to prove global stability, we need to show that $\bar V_1=0$. 

Recall that $V(x)=\sum_{i=1}^{\ab\Ss} \ab{\dot x_i}$, and $\sigma_i(t)=\sgn(x_i(t))$. Then, let the time-dependent sets $\Sigma_+,\Sigma_-,\bar\Sigma_+,\bar\Sigma_- \subset \{1,..,\ab\Ss\}$ be defined as: $ \Sigma_+(t)=\{i\vert \sigma_i(t)> 0\}$, $ \Sigma_-(t)=\{i\vert \sigma_i(t) < 0\}$, $\bar\Sigma_+(t)=\{i\vert \sigma_i(t)\ge 0\}$, and $ \bar\Sigma_-(t)=\{i\vert \sigma_i(t)\le 0\}$.

 Since $V(x)=\sum_{i=1}^{\ab\Ss}\ab{\dot x_i}$, we can write: 
\begin{equation}  \label{e.las1} V(x(t)) = \sum_{i\in \Sigma_+(t)} \dot x_i(t) - \sum_{i\in \Sigma_-(t)} \dot x_i(t) \equiv \bar V_1.  \end{equation}

If either one of the sets $\Sigma_+(t), \Sigma_-(t)$ is empty for some $t$, then this implies that $\dot x(t)\equiv 0$, and hence $\bar V_1=0$ which proves the statement. Therefore, we assume that both $\Sigma_+(t), \Sigma_-(t)$ are non-empty for all $t$.

Similar to the proof of Theorem \ref{mainth}, let $x_i^+$ denote the concentration of the dimer that corresponds to $X_i$, and let $x_j^-$ be the the concentration of the dimer that corresponds to $X_j$.   Denote $\dot R_{ij}(x_i,x_i^+,x_j^-):= \frac{\partial R_{ij}}{\partial x_i}\dot x_i + \frac{\partial R_{ij}}{\partial x_i^+}\dot x_i^+ + \frac{\partial R_{ij}}{\partial x_j^-}\dot x_j^-, \dot R_{i}(x_i,x_i^+):= \frac{\partial R_{i}}{\partial x_i}\dot x_i+\frac{\partial R_{i}}{\partial x_i^+}\dot x_i^+.$ Using the argument in the Theorem \ref{mainth}, a term of the form $R_{ij}(x)$ appears in $V$ with a positive coefficient only if $\sigma_i\ge 0$, $\sigma_j\le 0$, and $\sigma_i \sigma_j\ne0$. Similarly, $R_{ij}(x)$ appears in $V$ with a negative coefficient only if $\sigma_i\le 0$, $\sigma_j\ge 0$,  and $\sigma_i \sigma_j\ne0$. Hence, we can write:
\[\dot V= \sum_{(i,j)\in \bar \Sigma_- \times \bar\Sigma_+} \left( \rho_{ij} \dot R_{ij}(x) + \dot R_i(x)  \right) - \sum_{(i,j)\in \bar \Sigma_+ \times \bar\Sigma_-}\left( \rho_{ij} \dot R_{ij}(x) - \dot R_i(x) \right)\equiv 0, \]
for some $\rho_{ij}>0$. Note that the dependence on $t$ in the equation above has been dropped for notational brevity. 

As in the proof of Theorem \ref{mainth}, each term is nonpositive. Hence, $\dot V\equiv 0$ implies that \textit{each term} is identical to zero. We make several conclusions from the last statement:

First, $\forall i \in \Sigma_+ \cup \Sigma_-, \dot R_i(x(t))\equiv 0$. Furthermore, by definition,  $\dot R_i(x(t))=0$ for $i \not\in \Sigma_+ \cup \Sigma_i$. Therefore, we get that $\forall i, \dot R_i(x(t))\equiv 0$. Hence, 
$\sum_{i=1}^{\ab\Ss} \dot x_i(t) = \sum_{i=1}^{\ab\Ss} \sum_{j\ne i} (\dot R _{ji}(x(t))-\dot R_{ij}(x(t))) \equiv  0.$
Therefore, we get $ \sum_{i=1}^{\ab\Ss} \dot x_i(t) \equiv \bar V_2$ for some constant $\bar V_2$. Hence, 
\begin{equation}\label{e.las2}
	 V(x(t)) = \sum_{i\in \Sigma_+(t)} \dot x_i(t)+ \sum_{i\in \Sigma_-(t)} \dot x_i(t) \equiv \bar V_2.
	\end{equation}  

Second,  fix $i \in \Sigma_+(t)$. Then, for all $j\in\Sigma_-(t)$, we have $\dot R_{ij}(x(t)) \equiv 0$. Hence, we claim the following: If $i \in \Sigma_i^+(t')$ for some $t'>0$, then  $i \in \Sigma_i^+(t')$ for all $t>t'$. To show this, we can write 
\[\ddot x_i(t) = -F_i(\dot x_i) + \overbrace{ \sum_{\substack{i \in \Sigma_+ \\ j\ne i}} \dot R_{ji}(x)}^{>0} + \overbrace{\sum_{j\in \Sigma_-} \dot{R}_{ji}(x)}^{=0}, \]
where $F_i$  lumps all the terms that corresponds to reactions for which $  X_i$ is a reactant. It can be noted immediately that $F_i(\dot x_i)>0$ is positive when $\dot x_i>0$, and $F_i(0)=0$. The second term is positive since $R_{ji}$ is monotone and $\dot x_j>0$. The third term is zero as we have shown earlier. Therefore,  $\left . \ddot x_i\right\vert _{\dot x_i=0} >0$. Hence,  $\dot x_i(t)>0$ for all $t>t'$, i.e. $i \in \Sigma_+(t)$ for all $t\ge t'$ as claimed.

Similarly, if $i \in \Sigma_i^-(t')$ for some $t'>0$, then  $i \in \Sigma_i^-(t')$ for all $t>t'$.  Since $\Sigma_+,\Sigma_- \subset \{1,..,\ab\Ss\}$, there exists $T>0$ and constant sets  $\Sigma_+^*, \Sigma_-^*$ with $\Sigma_+(t)=\Sigma_+^*,\Sigma_-(t)=\Sigma_-^*$ for all $t\ge T$.
Therefore, we can write \eqref{e.las1},\eqref{e.las2} as: 
\[\sum_{i\in \Sigma_+^*} \dot x_i(t)- \sum_{i\in \Sigma_-^*} \dot x_i(t) \equiv \bar V_1, ~\sum_{i\in \Sigma_+^*} \dot x_i(t)+ \sum_{i\in \Sigma_-^*} \dot x_i(t) \equiv \bar V_2. \]
Hence, $\sum_{i\in\Sigma_+^*} \dot x_i(t) \equiv \tfrac 12 (\bar V_1+\bar V_2)$, and $\sum_{i\in\Sigma_-^*} \dot x_i(t) \equiv \tfrac 12 (\bar V_2-\bar V_1)$. By integrating the last two equations we note that  we can only have $\bar V_1=\bar V_2=0$ since $x(t)$ has been assumed to be bounded.
 
\strut \hfill $\square$

\subsubsection*{Proof of Theorem \ref{th.linearmod_J}}

Before proving the result, we need few preliminary lemmas. 

\begin{lemma}\label{lem.CB}
	Let $\Lambda=-\Gamma V$, where $\Gamma \in \mathbb R^{n \times \nu}$, $V \in \mathbb R^{\nu \times R}$. Let $r=\rank(\Lambda)$. Then,
	\begin{enumerate}
		\item Let $I$ be an arbitrary subset of $\{1,..,n\}$ with $\vert I\vert=r$. The corresponding principal minor can be written as:
		\begin{equation}\label{e.minor}\det_I(
		\Lambda)=\mbox{det}_I (-\Gamma V) = \sum_{\substack{J \subset \{1,..,\nu\}\\ \vert J\vert=r}} \det(-\Gamma_{IJ})\det(V_{JI}),\end{equation}
		where $\Gamma_{IJ},V_{JI}$ denote the submatrices of $\Gamma,V$ with the row and column indices specified in $I$ and $J$, respectively.
		\item The essential determinant of $\Lambda$ can be written as:
		\begin{equation}\label{essdet}
			\det_{ess}(\Lambda) = \sum_{\substack{I \subset \{1,..,n\}\\ \vert I\vert=r}}\sum_{\substack{J \subset \{1,..,\nu\}\\ \vert J\vert=r}} \det(-\Gamma_{IJ})\det(V_{JI}).
		\end{equation}
	\end{enumerate}
\end{lemma}
\begin{proof}
	It follows immediately from the Cauchy-Binet's formula \cite{horn_and_johnson} and the definition of the essential determinant.
\end{proof}
 
 Before stating the next lemma, we need some notation. Recall that  $\mathcal K_\Ns $ is the set of all possible Jacobian matrices $\partial R/\partial x$ defined on $\mathbb R_+^n$.  Hence,  any $V \in \mathcal K_\Ns $ can have $s$ nonzero entries which is equal to the number of reactant-reaction pairs.   Let $s$ be the number of nonzero entries of $V$, we list them as $v_1,..,v_s > 0$. 
Next, we  show that each term in the expansion \eqref{e.minor} is nonnegative. In other words, 
 \begin{lemma}\label{lem.terms}
	Let $\Ns$ be a  network that  is robustly $P_0$ and non-degenerate. Let the stoichiometry matrix be  $\Gamma \in \mathbb R^{n\times \nu}$, and let  $r=\rank(\Gamma)$.   Then,  $\forall I\subset \{1,..,n\}, \forall J\subset \{1,..,\nu\}$ with $\vert I\vert =\vert J\vert =r$ we have $ \det(-\Gamma_{IJ})\det(V_{JI}) \ge 0$ for all $V \in \mathcal K_\Ns$. 
\end{lemma}
\begin{proof}  As noted before, $v_1,..,v_s$ are the nonzero entries of  $V \in \mathcal K_\Ns$. Hence, we can think of $\det(V_{JI})$ as a polynomial in $v_1,..,v_s$. For the sake of contradiction, assume that there exists $I^*,J^*$ with $\ab I^*=\ab J^*=r$ such that $ \det(-\Gamma_{I^*J^*})\det(V_{J^*I^*})<0$ for some $V\in\mathcal K_\Ns$.
		 Hence, $\det(V_{J^*I^*})$ must have a monomial term $m^*=\prod_\ell v_\ell$ such that $\det(-\Gamma_{I^*J^*}) m^* <0$. 
		 
		 The   corresponding principal minor can be written using \eqref{e.minor} as $$\det_{I^*}(-\Gamma V) = \sum_{J, \ab J=r} \det(-\Gamma_{I^*J})\det(V_{JI^*}).$$
	 All the entries in $V$ that do not appear in $m^*$ can be set to be arbitrarily small. Since the determinant is  homogeneous in the entries of the corresponding matrix, all the terms other than $\det(-\Gamma_{I^*J^*}) m^* <0$ will become arbitrarily small in the expansion \eqref{e.minor}  which makes the principal minor $ \det_{I^*}(-\Gamma V)<0$ for some $V \in \mathcal K_\Ns$. However, 	  $\Ns$ is robustly $P_0$ which means that all principal minors of $-\G V$  are non-negative for any $V \in \mathcal K_\Ns$ which is a contradiction.
 \end{proof}

 We are ready now to prove the statement of the theorem.   Let $\tl\Ns=(\tl \Ss,\tl \Rs)$ be the modified network, and let $\tl\Gamma$ be its stoichiometry matrix.  Furthermore, since $\Ns$ is robustly non-degenerate, and using Corollary \ref{cor.nondeg} and Lemma \ref{lem.terms}, $\exists I\subset \{1,..,n\}, \exists J\subset \{1,..,\nu\}$ with  $\ab I=r$,$\ab J=r$ such that   $\det(-\Gamma_{IJ} )\det(V_{JI})>0$.
 
 The proof is divided by the graph modification under consideration.
\begin{enumerate}
	\item \emph{Reversal:}  Recall $\tl\Ss=\Ss$ and $\tl\Rs=\Rs \cup \{\R_{-j}\}$ for some $j$. Let $\Gamma_j$ denote the $j$th column of $\Gamma$. Then, $\tl\Gamma=[\Gamma,-\Gamma_j]$. Hence, $\rank(\tl\Gamma)=\rank(\Gamma)=r$. Using Corollary \ref{cor.nondeg}, we need show the existence of a positive $r\times r$ principal minor.    We consider the minor corresponding to $I$ for $\tl\Ns$. It can be seen immediately that is also positive since addition of a reverse of a reaction can only add new non-negative terms (by Lemma \ref{lem.terms})  to the expansion \eqref{e.minor}.
	
	\item \emph{Adding an intermediate:} Recall that the network $\Ns$ is modified by replacing a reaction of the form   $\R_j=\sum_{i}\alpha_{ij}X_i \to \sum_i \beta_{ij} X_i$ by two reactions $\tl\R_j:= (\sum_{i}\alpha_{ij}X_i \to X_{n+1})$, and $(\tl \R_{\nu+1}:=  X_{n+1} \to \sum_i \beta_{ij} X_i) $. W.l.o.g, assume that $\R_j$ is the last reaction, i.e., $j=\nu$. Let $\Gamma_\nu$ be the last column of $\Gamma$, hence we write $\Gamma=[\hat\Gamma,\Gamma_\nu]$. Define the two vectors $\Gamma_\nu^+,\Gamma_\nu^-$ entry-wise as $[\Gamma_\nu^+]_i=\max\{[\Gamma_\nu]_i,0\},  [\Gamma_\nu^-]_i=\min\{[\Gamma_\nu]_i,0\}$, respectively. In other words, $\Gamma_\nu^+$ contains the positive entries of $\Gamma_\nu^+$ while $\Gamma_\nu^-$ contains the negative entries of $\Gamma_\nu$. Hence, by construction, $\Gamma_{\nu}=\Gamma_\nu^++\Gamma_\nu^-$ Then, it can be seen that:
	\[ \tl\Gamma=\begin{bmatrix}\hat\Gamma & \tl\Gamma_\nu^- & \Gamma_\nu^+ \\ 0 & 1 & -1 \end{bmatrix}. \]
	Let $\tl I:=I \cup \{n+1\}, \tl J=J\cup \{\nu+1\}$. 	We will be computing $\det(-\tilde\Gamma_{\tl I\tl J})$ and $\det(\tl V_{\tl I\tl J})$. We start with the latter. The new reaction has only one reactant which is $  X_{n+1}$. Hence, we can write:  
	\begin{equation}\label{V1}\det(\tl V_{\tl J\tl I})=\det \left(\begin{bmatrix} V_{JI} & 0 \\ 0 & 1 \end{bmatrix} \right) = \det(V_{JI}). \end{equation}
		Next, we consider $\det(-\tilde\Gamma_{\tl I\tl J})$. We study two cases: $\nu \in J$ and $\nu \not\in J$. 
	
	\noindent\underline{Case 1: $\nu \in J$:} Let $J^*=J/\{\nu\}$. We get:
 	\begin{align}\nonumber\det(-\Gamma_{\tl I \tl J})&=\det\left(-\begin{bmatrix} \hat \Gamma_{IJ^*} & \tilde \Gamma_{\nu,I}^- & \tilde \Gamma_{\nu,I}^+ \\ 0 & 1 & -1 \end{bmatrix}\right) \mathop=^{(\heartsuit)} \det\left(-\begin{bmatrix} \hat \Gamma_{IJ^*} & \tilde \Gamma_{\nu,I} & \tilde \Gamma_{\nu,I}^+ \\ 0 & 0 & -1 \end{bmatrix}\right) \\ \label{G1} & = \det\left(-\begin{bmatrix}  \Gamma_{IJ} &  \tilde \Gamma_{\nu,I}^+ \\ 0 & -1 \end{bmatrix}\right)= \det(-\Gamma_{IJ}),  \end{align}
 	where $(\heartsuit)$ follows by the fact that the determinant is invariant under the addition of the last two columns. \\
 	\noindent\underline{Case 2: $\nu \not\in J$:} We get:
 	\begin{align}\label{G2} \det(-\Gamma_{\tl I \tl J})&=\det\left(-\begin{bmatrix}  \Gamma_{IJ} &  \tilde \Gamma_{\nu,I}^+ \\ 0  & -1 \end{bmatrix}\right) = \det(-\Gamma_{IJ}),  \end{align}

Hence, using \eqref{V1},\eqref{G1},\eqref{G2}, we get that $\det(-\tilde \G_{\tl I\tl J} \tilde V_{\tl J\tl I})=\det(- \G_{ I J} V_{ J I}) > 0$. Finally, since $\tl\Ns$  admits an RLF, then Lemma \ref{lem.terms} and Corollary \ref{cor.nondeg} imply that the $(r+1)\times(r+1)$ principal minor corresponding to $\tl I$ for modified network $\tl\Ns$ is positive, hence $\tl\Ns$ is robustly non-degenerate. 

\item \emph{External regulation of a species:} $\tl\Ss=  \Ss$, and $\exists X_i \in \Ss$ such that $\tl\Rs=\Rs \cup \{ X_i \leftrightharpoons \emptyset \}$. W.l.o.g, assume that $i=n$. Then we study two cases: $\rank(\tl\Gamma)=\rank \Gamma$ and  $\rank(\tl\Gamma)=\rank \Gamma+1$. In first case, we let $\tl I=I, \tl J=J$. Hence, $\det(-\tl\Gamma_{\tl I \tl J})\det(\tl V_{\tl J \tl I})=\det(- \Gamma_{I  J})\det( V_{  JI})>0$. Therefore, using Lemma \ref{lem.terms} and Corollary \ref{cor.nondeg}, $\tl\Ns$ is robustly non-degenerate. 

We now study the case in which  $\rank\tl\Gamma=1+\rank\Gamma=1+r$. We will first claim that it must be possible to choose $I$ such that $n\not\in I$. As a proof, consider the contrary. Then, this means that  $\forall I\subset\{1,..,n\}$ that satisfies $\ab I=r$ and $n\not \in I$, we have $\det_I(-\Gamma V)=0$. Using Lemma \ref{l.essdet}, this means that removing   $\gamma_n^T$ (the $n$th row of $\Gamma$), i.e., removing $X_n$ from $\Ns$, will cause the rank of $\Gamma$ to drop from $r$ to $r-1$. Hence, this means that $\gamma_n^T$ is linearly independent from the other rows of $\Gamma$. However, adding the reactions $ \{ X_i \leftrightharpoons \emptyset \}$ will only modify the $n$th row in $\Gamma$. Since $\gamma_n$ is already independent of the other rows of $\Gamma$, the rank cannot increase, which is a contradiction. 

Therefore, when $\rank\tl\Gamma=r+1$, we let $I$ be chosen such that $n\not\in I$. Hence,  let $\tl I=I \cup \{n\}, \tl J=J\cup\{\nu+1\}$, where $\tl\R_{\nu+1}:=   X_n \to\emptyset$. Therefore, we can write, 
	\begin{align}\label{G3} \det(-\tl\Gamma_{\tl I \tl J})&=\det\left(-\begin{bmatrix}  \Gamma_{IJ} &  0 \\ 0  & -1 \end{bmatrix}\right) = \det(-\Gamma_{IJ}),  \end{align}
Similarly, $\det(\tl V_{\tl J \tl I})=\det( V_{ J I})$ using the same argument as in \eqref{V1}. Therefore, $\det(-\tl\Gamma_{\tl I \tl J})\det(\tl V_{\tl J \tl I})=\det(- \Gamma_{I  J})\det( V_{  JI})>0$. Finally, using Lemma \ref{lem.terms} and Corollary \ref{cor.nondeg}, $\tl\Ns$ is robustly non-degenerate.  
		\item \emph{Conserved Regulation of a species:} We can consider this case as a sequence of two modifications. First, let $\hat\Ss=  \Ss$, and $\hat\Rs=\Rs \cup \{ X_i \leftrightharpoons \emptyset \}$. W.l.o.g, assume that $i=n$. Then, from the previous case it follows that $\hat \Ns=(\hat\Ss,\hat\Rs)$ is robustly non-degenerate. Let $\hat\Gamma$ be the corresponding stoichiometry matrix. Hence, using Lemma \ref{lem.terms} and Corollary \ref{cor.nondeg}, there exist sets $\hat I, \hat J$ with $\ab {\hat I}=\ab {\hat J}=\rank{\hat\G}$ such that $\det(-\hat\G_{\hat I \hat J})\det(\hat V_{\hat J \hat I})>0.$
		
		 Next, we define $\tl\Ns$ as follows: $\tl\Ss=\hat\Ss\cup\{  X_{n+1}\} $ and $\tl\Rs$ is defined as follows: both $\tl\Rs,\hat\Rs$ have the same reactions except for $X_n \leftrightharpoons \emptyset$ which is replaced by $X_n \leftrightharpoons X_{n+1}$.  Now consider two cases:   $\rank(\tl\Gamma)=\rank({\hat\Gamma})$ and  $\rank(\tl\Gamma)=1+\rank{(\hat\Gamma)}$. In first case, using the same argument as in the case of external regulation, $\tl\Ns$ is robustly non-degenerate.
		
		We now study the case in which  $\rank\tl\Gamma=1+\rank{\hat\G}$. Using a similar argument to the case of external regulation, we can choose $\hat J$ such that $\nu,\nu+1 \not\in \hat J$, where $\R_{\nu},\R_{\nu+1}$ are the reactions $  X_{n+1}\to   X_n,   X_n \to   X_{n+1}$, respectively.  Hence,  let $\tl I=\hat I \cup \{n+1\},  \tl J=\hat J\cup\{\nu\}$. 
		Hence, similar to \eqref{G3} and \eqref{V1} we can show that $\det(-\tl\G_{\tl I \tl J})= \det(-\tl\G_{\hat I \hat J})$, $\det(\hat V_{ \tl J \tl I})= \det(\hat V_{\hat J \hat I })$.  Therefore, $\det(-\tl\Gamma_{\tl I \tl J})\det(\tl V_{\tl J \tl I}))>0$. Finally, using Lemma \ref{lem.terms} and Corollary \ref{cor.nondeg}, $\tl\Ns$ is robustly non-degenerate.

	\item \emph{Adding a catalyst:} Let $\tl\Ss=\Ss\cup\{  X_i^-\}$, and $\tl\Rs$ is defined as in Definition \ref{def.mod}, Item 6. It can be seen that this implies that $\tl \G=[\G^T,-\gamma_i]^T$, where $\gamma_i$ is the $i$th row of $\G$. Therefore $\rank\tl\Gamma=\rank\Gamma$. Let $\tl I=I, \tl J=J$.  Since all the reactions in $\tl\Ns$ are extensions of the corresponding reactions in $\Ns$, then the positive term $\det(-\Gamma_{IJ})\det(V_{JI})$ is present in the expansion of $\det_{\tl I}(-\tl \Gamma \tl V)$. Therefore,  using Lemma \ref{lem.terms} we get that $\det_{\tl I}(-\tl \Gamma \tl V)>0$. Hence, using Corollary \ref{cor.nondeg}, $\tl \Ns$ is robustly non-degenerate.
	\item \emph{Adding a dimer:} Let $\tl\Ss=\Ss\cup\{  X_i^+\}$, and $\tl\Rs$ is defined as in Definition \ref{def.mod}, item 7. It can be seen that this implies that $\tl \G=[\G^T,\gamma_i]^T$, where $\gamma_i$ is the $i$th row of $\G$. Therefore $\rank\tl\Gamma=\rank\Gamma$. Using the same argument as in the previous case we get that $\tl \Ns$ is robustly non-degenerate.

 \end{enumerate}
  \strut \hfill $\square$
 
 {
 \subsubsection*{Proof of Theorem \ref{th_cs1} }
 Let $\Ns=(\Ss,\Rs)$, and let $\tl\Ns=(\tl\Ss,\tl\Rs)$ be its elementary modification. For a given reaction $\R_j$, let $\mathcal I(\R_j)\subset \Ss$ denotes its reactants, while $\mathcal O(\R_j) \subset \Ss$ denotes its products. The statement of the theorem is equivalent to proving that the absence of critical siphons for $\Ns$ implies the same for $\tl\Ns$.  Pick any $P \subset \Ss$. By assumption, $P$ is not a critical siphon for $\Ns$. Hence, $P$ is either not a siphon, or it is a trivial siphon. For the first case, using the definition of a siphon, $P$ is not a siphon if and only if  the following statement $(\clubsuit)$ holds : \begin{center}
 	\begin{minipage}{0.8\textwidth}$(\clubsuit)$ $\exists X_k \in P, \R_k \in \Rs$ such that $X_k$ is a product of $\R_k$ (i.e, $X_k\in \mathcal O(\R_k)$), and $\mathcal I(\R_k)\cap P =\emptyset$.\end{minipage}
 \end{center}  For the second case, if $P$ is a trivial siphon,  we assume, w.l.o.g, that it is minimal, i.e., $P$ coincides exactly with the support of a single conservation law. %
   We are ready now to consider the following cases:
 \begin{enumerate}	\item \underline{Reversal:} By definition, $\tl\Ss=\Ss$,  and for some $j\in \{1,..,\nu\}$ we have $\tl\Rs=\Rs\cup\{\R_{-j}\}$.  If $P$ is not a siphon for $\Ns$, and since $\tl\Ss=\Ss$ (i.e, no new species added), then the statement $(\clubsuit)$ holds also for $\tl\Ns$. Hence, $P$ is not a siphon for $\tl\Ns$.  If $P$ is a trivial siphon for $\Ns$, then it is also a trivial siphon for $\tl\Ns$ since addition of a reverse of a reaction does not change the conservation laws of a network. In summary, $P$ is not a critical siphon for $\tl\Ns$. Since $\tl \Ss=\Ss$,  $\tl\Ns$ lacks critical siphons.
 	\item \underline{External regulation:}  $\tl\Ss=  \Ss$, and $\exists X_i \in \Ss$ such that $\tl\Rs=\Rs \cup \{ X_i \leftrightharpoons \emptyset \}$. Pick any $P \subset \Ss$.  If $P$ is not a siphon, then the same argument used for the previous modification shows $P$ is not a siphon for $\tl\Ns$. If $P$ is a (minimal) trivial siphon for $\Ns$, then either:   (A) $X_i \not\in P$ which means that $P$ is a trivial siphon for $\tl\Ns$, or (B) $X_i \in P$ which means that $P$   no longer contains the support of a conservation law  for $\tl\Ns$ since $X_i$ has an inflow and is no longer conserved.  In summary, $P$ is not a critical siphon for $\tl\Ns$. Since $\tl \Ss=\Ss$,    $\tl\Ns$ lacks critical siphons.
 	
 	\item \underline{Conserved regulation:} $\tl\Ss=  \Ss \cup \{X^*\}$, and $\exists X_i \in \Ss$ such that $\tl\Rs=\Rs \cup \{ X_i \leftrightharpoons X^* \}$.  If $P$ is not a siphon, then $P$ is not a siphon for the modified network $\tl\Ns$ because the statement $(\clubsuit)$ continues to hold. If we define $\tl P:=P \cup \{X^*\}$, then $\tl P$ is not a siphon since $\mathcal{I}(\R_k)\cap \tl P =\emptyset$, i.e. $(\clubsuit)$ holds. If $P$ is a (minimal) trivial siphon for $\Ns$, then either:   (case A) $X_i \not\in P$ which means that $P$ is a trivial siphon for $\tl\Ns$, or, (case B) $X_i \in P$ which means that $P$ no longer contains the support of a conservation law for $\tl\Ns$ and  $P$ is no longer a siphon for $\tl\Ns$.  Instead, $P\cup\{X^*\}$ contains the support of a conservation law, and hence $P\cup\{X^*\}$ is a trivial siphon for $\tl\Ns$.  In summary, neither $P$ nor $P\cup\{X^*\}$ are critical siphons for $\tl\Ns$.  Since all  subsets of $\tl\Ns$ can be represented as $P$ or $P \cup \{X^*\}$ for some $P\subset \Ss$,  $\tl\Ns$ lacks critical siphons.
 	
 	\item \underline{Adding an intermediate:}  Recall that the network $\Ns$ is modified by replacing a reaction of the form   $\R_j=\sum_{i}\alpha_{ij}X_i \to \sum_i \beta_{ij} X_i$ by two reactions $\tl\R_j:= (\sum_{i}\alpha_{ij}X_i \to X_{n+1})$, and $(\tl \R_{\nu+1}:=  X_{n+1} \to \sum_i \beta_{ij} X_i) $.  If $P$ is not a siphon for $\Ns$, then the statement $(\clubsuit)$ continues to hold for $\tl\Ns$, i.e., $P$ is not a siphon for $\Ns$.  Next, let $\tl P=P\cup\{X^*\}$. For the sake of contradiction, assume that $\tl P$ is a siphon, this is only possible if $X^*\in \mathcal I(\R_k)$ (where $\R_k$ is defined in the statement $(\clubsuit)$). But using our construction, this means that $\mathcal \R_k=\R_j$. Furthermore, since $\tl P$ is a siphon and it contains $X^*$, one of the reactants of $\tl \R_j$ is in $\tl P$. This means that one of the reactants of $\R_j(=\R_k)$ is in $P$. But, this contradicts the statement $(\clubsuit)$. Hence, $\tl P$ is not a siphon for $\tl\Ns$.  %
 	
 	If $P$ is a (minimal) trivial siphon for $\Ns$, then it contains the support of a conservation law $d \in \mathbb R^n_{\ge 0}$. Since $P$ is the support of $d$, denote $s:=\vert  P\vert$, and recall that $\gamma_{ij}$ denotes the $(i,j)$th entry of $\Gamma$. W.l.o.g, assume that $\Ss$ is indexed such that the first $s$ elements coincide with the elements of $P$. This also implies that  $d_1,..,d_s >0$. Hence,   $\forall j \in \{1,..,\nu\},~\sum_{i=1}^s d_i \gamma_{i j} =0$. Next, we consider few cases: (Case A) $\mathcal I(\R_j) \cap P =\emptyset, \mathcal O(\R_j) \cap  P =\emptyset$. We can see that the addition of an intermediate does not change the conservation law, therefore $P$ is a trivial siphon for $\tl\Ns$, while $P \cup \{X^*\}$ is not a siphon since we assumed that $\mathcal I(\R_j)\cap P =\emptyset$. (Case B) $\mathcal I(\R_j) \cap  P =\emptyset, \mathcal O(\R_j) \cap P \ne \emptyset$, or $\mathcal I(\R_j) \cap  P \ne \emptyset, \mathcal O(\R_j) \cap  P =\emptyset$. By AS2, this means either $\sum_{1\le i \le s, \gamma_{ij}>0} d_i \gamma_{i j} =0$ or $\sum_{1\le i \le s, \gamma_{ij}<0} d_i \gamma_{i j} =0$, respectively. Either case contradicts   $d_1,..,d_s>0$. (Case C) $\mathcal I(\R_j) \cap P \ne\emptyset, \mathcal O(\R_j) \cap  P \ne \emptyset.$ Using AS2, we get $\sum_{1\le i \le s, \gamma_{ij}<0} d_i \gamma_{i j}=\sum_{1\le i \le s, \gamma_{ij}>0} d_i \gamma_{i j}:=\xi>0$. Therefore $\tl P=P\cup \{X^*\}$ is a trivial siphon for $\tl\Ns$ with the conservation law $\tilde d=[d^T \xi]^T$, while $P$ is not a siphon for $\tl\Ns$.  In summary, neither $P$ nor $P\cup\{X^*\}$ are critical siphons for $\tl\Ns$.  Since all  subsets of $\tl\Ns$ can be represented as $P$ or $P \cup \{X^*\}$ for some $P\subset \Ss$,  $\tl\Ns$ lacks critical siphons.

 	\item \underline{Adding a dimer:}  Let $\tl\Ss=\Ss\cup\{  X_i^+\}$, and $\tl\Rs$ is defined as in Definition \ref{def.mod}, item 7. If $P$ is not a siphon, this means that the statement $(\clubsuit)$ holds. Since $X_i^+$ shares the same input and output reactions with $X_i$,  the statement $(\clubsuit)$ continues to hold. Hence, $P$ is not a siphon for $\tl\Ns$. If $P$ is a trivial siphon for $\Ns$ then it will be a trivial siphon for $\tl\Ns$ since   adding a dimer preserves the existing conservation laws of $\Ns$.
 	
 	Now, let us consider a set of the form $\tl P:=P\cup\{X_i^*\},P\subset \Ss$. For the sake of contradiction, assume that $\tl P$ is a critical siphon for $\tl \Ns$.  Let $\hat P:= (\tl P/\{X_i^+\})\cup\{X_i\}$. Since $X_i$ have the same reactants and products as $X_i^+$ by construction, then $\hat P$ is a siphon for $\Ns$ and it does not contain the support of a conservation law . Hence, $\hat P$ is a critical siphon for $\Ns$ which contradicts our assumption.  Since all  subsets of $\tl\Ns$ can be represented as $P$ or $P \cup \{X^*\}$ for some $P\subset \Ss$,  $\tl\Ns$ lacks critical siphons.
\end{enumerate}
 
\subsubsection*{Proof of Theorem \ref{th_cs2}-1}
 Part 1 of Theorem \ref{th_cs2} follows from Theorem \ref{th_cs1} except for the case of adding a catalyst which is proved next. Assume that $\tl\Ns$ is a modification of a linear network $\Ns$ by adding a catalyst as described in Definition \ref{def.mod}. Let $P \subset \Ss$. So $P$ is either not a siphon or a trivial siphon. Assume first that $P$ is not a siphon, hence the statement $(\clubsuit)$ in proof of Theorem \ref{th_cs2} holds. Since the products and reactants of reactions in $\tl\Ns$ contain their counterparts in $\Ns$ then $P$ is not siphon for $\tl\Ns$. Consider $\tl P=P\cup\{X_i^-\}$. We consider two cases: (case A) $X_i \not\in \tl P$. Note both $\Ns$ and $\tl\Ns$ satisfy AS1. Therefore, $\tl\Rs$ must contain at least one reaction $\tl\R_j$ with $X_i^- \in \mathcal O(\tl\R_j)$.  In order for $\tl P$ to be a siphon, there must exist $X^* \in \mathcal I(\tl\R_j) \cap \tl P$. But since $\Ns$ is linear, $\R_j$ has at most one reactant and it must be $X_i$ (since by construction, input reactions of $X_i^-$ are output reactions of $X_i$). This also implies that $\tl\R_j$ has at most one reactant, and hence $X_i=X^*\in \tl  P$ which is a contradiction.   (case B) $X_i \in \tl P$, hence, $\tl P$ contains the support of the conservation law $X_i + X_i^-=\mbox{constant}.$ Therefore, $\tl P$ is a trivial siphon, and it is not critical. 
 
 Second, let us assume that $P$ is a trivial siphon for $\Ns$, then $P$ and $P\cup\{X_i^-\}$ are trivial siphons for $\tl\Ns$ since   adding a catalyst preserves the existing conservation laws of $\Ns$.  Since all  subsets of $\tl\Ns$ can be represented as $P$ or $P \cup \{X^*\}$ for some $P\subset \Ss$,  $\tl\Ns$ lacks critical siphons. 
 
}
 \section*{Declarations}
 \noindent \emph{Funding:} This research has been funded by NSF grant 2052455. \\
\noindent \emph{Conflict of interest:} The author declares no conflict of interest. \\
\noindent \emph{Author Contributions:} This article has a single author who performed all related tasks. \\

\end{document}